\documentclass[preprint,12pt]{elsarticle}

\pdfoutput=1
\usepackage{amsmath,amsbsy,amssymb,color}
\usepackage{epsfig}
\usepackage{mathptmx}
\usepackage{latexsym}
\usepackage[unicode=true, colorlinks=true, linkcolor=blue]{hyperref}
\usepackage{url}
\usepackage{float}
\usepackage{natbib}
\newdefinition{rmk}{Remark}

\journal{journal to be defined}

\begin{document}
\begin{frontmatter}

\title{A modified DFC method based on full state and periodic gain for stabilization and chaos control}

\author{Ver\'onica E. Pastor \& Graciela A. Gonz\'alez}

\begin{abstract} {A novel delayed feedback control based on full state is proposed. The designed scheme combines the difference between  two delayed states  and a periodic control gain. System stabilization is achieved in any hyperbolic unstable equilibrium point. The procedure to build the control is systematic and, the set of the stabilizing control parameters is rigorously featured by analytical arguments. An ad-hoc strategy based on the proposed scheme is implemented for control of chaos by stabilizing equilibrium points embedded in chaotic attractors. Simulated examples illustrate this strategy performance.  }

\end{abstract}

\end{frontmatter}

\section{Introduction} \label{intro}

Since Pyragas proposed the delayed feedback control (DFC) method in 1992 (\cite{PYRAGAS1992421}), there has been increasing  interest in stabilizing unstable periodic orbits (UPO) via this resource. A distinct feature of  DFC and its modifications, is that it does not require preliminary knowledge about the desired UPO except its period. The UPO of the free system remains as a periodic orbit of a controlled system by a DFC method. These algorithms are noninvasive that is, the control perturbation vanishes whenever the system settles on the target orbit. 
The extended DFC (EDFC) was introduced in \cite{Socolar} for dealing with highly unstable orbits and it has resulted succesfull on some applications as reported in  \cite{Zheng-Yu-2022}. Besides, a dynamic EDFC has been presented in \cite{Zheng-Yu-2022} as an improvement over the classical EDFC for stabilizing multi-rotation UPOs by enlarging the stability regions. 
The theoretical analysis of the DFC method comes out difficult because the resulting dynamics evolves in an infinite dimensional phase space. The introduction of a periodic control gain, based on the act-and-wait concept, makes the stability study hing on a finite dimensional monodromy matrix. A great deal of works in this line are referred in \cite{PYRAGAS2018574} and there is a more recent paper (\cite{Liu2020}) which advances on numerical aspects for obtaining stability conditions of the closed loop system. As a further development of the EDFC, it is worth to mention \cite{Py-Py2019} where a state-dependent switch is appended to the feedback loop.

The DFC method has been reformulated for stabilizing a given dynamical system in one of its unstable equilibrium point. Namely it appears as an interesting tool for controlling chaos. It is designed as an amplified difference between the present state $x(t)$ and the delayed state $x(t-\tau)$ or, between a suitable observable of the state at times $t$ an $t-\tau$, so depending on the time delay $\tau$ and a control gain $K$. The extended version of DFC (EDFC) includes not only $\tau$ but also its integer multiples and its design is based on a geometric sum of feedback signals {(\cite{Chang1998}). Then, the multiple DFC (MDFC)  which depends on $n$ different delay $\tau_1$, $\tau_2$, \dots, $\tau_n$ and different gains $k_1$, $k_2$, \dots, $k_n$ was proposed in \cite{Ahlborn-Parlitz2005} as an alternative that offers major stability and flexibility over its antecedents.

As the controlled system by DFC, EDFC or MDFC involves a nonlinear delayed-differential equation, its stability analysis comes out rather intricate. Adequate values of $K$ and $\tau$  may be obtained by numerical computation but an analytical determination of them  results a critical point. In some works, the stabilizing parameters are obtained by solving some linear matrix inequalities. It may be mentioned  \cite{GUAN2006} which deals with the DFC where the control gain $K$ is an $n \times n$ matrix, being $n$  the phase space dimension, namely it is a delayed full state feedback control. An analytical procedure based on inequalities is also presented in \cite{Lu-Li2009} to obtain the set of stability parameters of a particular case of the MDFC.

On the other hand, if the linearization of the system in the equilibrium has either a zero eigenvalue or an odd number of real positive eigenvalues, there are no control parameters values  such that stabilization is achieved by any of the referred methods. This restriction is known as the odd number limitation (ONL) ( \cite{15H}, \cite{KUZNETSOV2015706}, \cite{Shumafov}).  
The introduction of a periodic  control gain allows to study the dynamics of the DFC by means of a discrete-time system. A scheme including  a periodic  control gain appears in \cite{Konishi2011}, based on the act-and.wait-concept. It is fully developed for the unstable focus case but  the ONL is not dealt with by this proposal. Later, a DFC algorithm based on a periodic  control  gain has been developed in \cite{Lenov2018} which overcomes the ONL but it only covers  the case of saddles which eigenvalues have unit geometric multiplicity. 
 
In \cite{VEP-GAG}, the classical difference of the DFC is replaced by the difference between the two delayed states $x(t-2\tau)$ and $x(t-\tau)$. To emphasize this modification, this scheme will be referred here as two-delayed feedback control (TDFC). It has been proved in \cite{VEP-GAG} that stabilization is achieved in the general nonlinear scalar case by TDFC with a periodic control gain. The periodic control law is simple and the adequate control parameters can be  computed systematically. 
 
This work deals with the extension of the results from \cite{VEP-GAG} to the n-dimensional case. The designed scheme combines a full state TDFC and a periodic control gain. By this proposal, system stabilization is achieved in any hyperbolic unstable equilibrium point (saddle, focus o saddle-focus) so in particular, ONL is overcome. Moreover, the procedure to build the control is systematic and, the set of the stabilizing control parameters is rigorously featured by analytical arguments.
Finally, the issue of controlling chaos is considered. An ad-hoc strategy based on the proposed scheme is introduced for stabilizing equlibrium points embedded in chaotic attractors. Simulated examples illustrate this strategy performance.

\section{Problem Statement}

Let us consider a system given by:\\
\begin{equation} \label{eq1}
    \dot{x}=f(x) 
\end{equation} 
where $x \in \mathbb{R}^n$ is the state vector and $f: \mathbb{R}^n \rightarrow \mathbb{R}^n$ is a piece-wise continuously differentiable function.\\
System (\ref{eq1}) has an unstable hyperbolic equilibrium point $x^*$. The objective is to stabilize the system in $x^*$ by applying additive control, this is,\\
\begin{equation} \label{eq2}
    \dot{x}=f(x)+u(t) 
\end{equation}
with $u(t)$ a time-delayed feedback control.\\

The control law is proposed as a generalization of the TDFC method  developed in \cite{VEP-GAG} to the n-dimensional case as follows, 

\begin{equation}\label{eq3}
    u(t)=K(t)\big(x(t-2\tau)-x(t-\tau)\big) \\
\end{equation}
where $\tau > 0$ and \\
\begin{equation}\label{eq3bis}
K(t)= \begin{cases}
0_{n \times n} \! & \mbox{if }  3k\tau \leqslant t < (3k+2)\tau \\
K \! & \mbox{if }  (3k+2)\tau \leqslant t < (3k+3)\tau \\
\end{cases} \ \ \ k \in \mathbb{N} \cup \{0\},
\end{equation}
being $0_{n \times n}$ the null matrix in $\mathbb{R}^{n \times n}$ and $K$ a constant matrix in $\mathbb{R}^{n \times n}$ which will be chosen adequately. 
 
Note that $x^*$ remains as an equilibrium point of (\ref{eq2}). The problem is finding $K$ and $\tau$ such that $x^*$ becomes an asymptotic stable equilibrium point of the controlled system (\ref{eq2})-(\ref{eq3}). \\
By putting $\delta x=x-x^*$, system (\ref{eq2})-(\ref{eq3}) results:
\begin{equation} \label{eq5}
    \delta \dot{x}=A \delta x+g(\delta x)+K(t)\big(\delta x(t-2\tau)-\delta x(t-\tau)\big)
\end{equation}
where $A=Df(x^*)$ is the jacobian matrix of $f$ in $x^*$ and $g$ is the function that comprises the nonlinear terms with $g(0)=0$.\\
In turn, there exists an invertible matrix $V$ such that\\
\begin{equation*}
        A=V^{-1}DV
\end{equation*}
being $D$ the Jordan matrix of $A$. Let us point out that as $x^*$ is hyperbolic, each eigenvalue of $A$ (or equivalently, of $D$) has real part different from zero. 

Introducing the coordinate change:\\
\begin{equation*} \label{eq6}
    y=V\delta x
\end{equation*}
and defining: 
\begin{equation}\label{eq7} \widetilde{K}=VKV^{-1} \end{equation} 
system (\ref{eq5}) is re-writen as:\\
\begin{equation*}
    \dot{y}=Dy+Vg(V^{-1}(y))+\widetilde{K}(t)\big(y(t-2\tau)-y(t-\tau)\big)
\end{equation*}
where \\
\begin{equation}\label{eq8bis}
\widetilde{K}(t)= \begin{cases}
0_{n \times n} \! & \mbox{if }  3k\tau \leqslant t < (3k+2)\tau \\
\widetilde{K} \! & \mbox{if }  (3k+2)\tau \leqslant t < (3k+3)\tau \\
\end{cases} \ \ \ k \in \mathbb{N} \cup \{0\},
\end{equation}

Working out as in the scalar case (see \cite{VEP-GAG} for details), it is valid to neglect nonlinear terms and to look for values $\tau$ and $\widetilde{K}$ that stabilize the system.
\begin{equation*}
    \dot{y}=Dy+\widetilde{K}(t)\big(y(t-2\tau)-y(t-\tau)\big)
\end{equation*}
in the origin. Then, by means of the relationship (\ref{eq7}), the suitable $\tau$ and $K$ are obtained.\\

Before developing this paper proposal, an outline on the stated results in \cite{VEP-GAG} for the scalar case will be re-visited. Consider the scalar controlled system defined as:\\
\begin{equation*}
    \dot{x}(t)=f(x(t))+ \epsilon(t)\big(x(t-2\tau)-x(t-\tau)\big)
\end{equation*}
with\\
\begin{equation} \label{eq7bis}    \epsilon (t)= 
\begin{cases}
0, & \mbox{if } 3k\tau\! \leqslant \!t <\! (3k\!+\!2)\tau \\
\epsilon, & \mbox{if } (3k\!+\!2)\tau\! \leqslant \!t<\!(3k\!+\!3)\tau
\end{cases} \ \ \ k \in \mathbb{N} \cup \{0\}.\\
\end{equation}

Let $x^*$ an equilibrium point of $f$ with $f'(x^*)=\lambda>0$. For each fixed $\zeta: \left|\zeta\right|<1$, the stabilizing $\tau$ and $\epsilon$ relates as:
\begin{equation}\label{pp2-1}
         \epsilon=\frac{e^{-\lambda\tau} (e^{3\lambda\tau}-\zeta)}{\tau (e^{\lambda\tau}-1)}
     \end{equation}


\section{The two-dimensional case}

The controlled system (\ref{eq5}) takes the form:\\
     \begin{equation*}\label{}
         \begin{pmatrix}
				\dot{x_1} \\ \dot{x_2} \end{pmatrix} = A\begin{pmatrix}
				x_1 \\ x_2 \end{pmatrix} +g\begin{pmatrix}
				x_1 \\ x_2 \end{pmatrix} + K(t)\begin{pmatrix}
				x_1(t-2\tau)-x_1(t-\tau) \\ x_2(t-2\tau)-x_2(t-\tau) \end{pmatrix}
     \end{equation*}
being $A$ the jacobian matrix of the system in its equilibrium point $\begin{pmatrix}
				x_1^* \\ x_2^* \end{pmatrix}$ and $K(t)$ defined as in (\ref{eq3bis}) for $n=2$.\\ 
    
The existence of solution is shown by well known analytic tools.
The characterization of the stability control parameters will be stated. Three cases will be analyzed, by considering the three possible Jordan forms of the matrix $A$:\\ 
\begin{equation*} 
(3.1) \; D= \ \ \begin{bmatrix} 
\lambda_1 & 0 \\
0 & \lambda_2 
\end{bmatrix} \ \  \  \
(3.2) \;  D= \ \ \begin{bmatrix} 
\lambda & 1 \\
0 & \lambda 
\end{bmatrix} \ \  \   \
(3.3) \; D= \ \ \begin{bmatrix} 
\mu & -\omega\\
\omega & \mu\\
\end{bmatrix} 
\end{equation*}

\vspace{1cm}

\textbf{\textsl{Case 3.1}}

\vspace{.25cm}
As $x^*$ is hyperbolic, both $\lambda_1$ and $\lambda_2$ are different from zero and as it is an unstable equilibrium point, at least one of the diagonal values are greater than zero. It is deduced straightfully from (\ref{pp2-1}) that the stabilization problem is solved by putting:

\begin{equation*}
\widetilde{K}= \begin{bmatrix} 
k_1 & 0 \\
0 & k_2 \\
\end{bmatrix} 
\end{equation*}

being \begin{equation*} k_1= \begin{cases} 0 \! & \mbox{if }  \lambda_1 <0\\
\epsilon_1 \! & \mbox{if }  \lambda_1 >0 \end{cases} 
\ \ \ \ \mbox{and} \ \ \ \    
k_2= \begin{cases}0 \! & \mbox{if } \lambda_2 <0\\
\epsilon_2 \! & \mbox{if } \lambda_2 >0 \end{cases} 
\end{equation*}

where \\
         \begin{equation}\label{eq10}
\epsilon_i=\displaystyle\frac{e^{-\lambda\tau} (e^{3\lambda\tau}-\zeta_i)}{\tau (e^{\lambda\tau}-1)} \ \ \ i=1,2.
         \end{equation}
with fixed $\zeta_i: \left|\zeta_i\right|<1$
         
\vspace{1cm}     

\textbf{\textsl{Case 3.2}}

\vspace{0.25cm}
From the instability of $x^*$, it results $\lambda >0$. It is easy to prove that stability is achieved by choosing\\
\begin{equation*} 
   \widetilde{K}=
   \begin{bmatrix} \epsilon &  0 \\
0 & \epsilon  \end{bmatrix}
\end{equation*}
where $\epsilon$ verifies (\ref{pp2-1}). 

As in the proof of Proposition 2 in \cite{VEP-GAG} the stabilization results lay on building a map $P$ such that\\
\begin{equation*} P\begin{pmatrix}
   x_1(3k) \\ x_2(3k) \end{pmatrix} = \begin{pmatrix} x_1(3k+3) \\ x_2(3k+3) \end{pmatrix},
\end{equation*} 
and it comes out: \\

\begin{equation} \label{eq11}
   P'(0)= \begin{bmatrix} 
      e^{3\lambda \tau}+\epsilon \tau e^{\lambda \tau}(1-e^{\lambda \tau})  & 3\tau e^{3 \tau \lambda}-2 \epsilon \tau ^2 e^{2 \lambda \tau}\\
      0  & e^{3\lambda \tau}+\epsilon \tau e^{\lambda \tau}(1-e^{\lambda \tau})
   \end{bmatrix} 
\end{equation}

By fixing a value $\zeta: |\zeta|<1$ and computing $\epsilon$ through the formula (\ref{pp2-1}), $\zeta$ is a double eigenvalue of $P'(0)$, yielding to the stabilization of the equilibrium point.\\


\textbf{\textsl{Case 3.3}}

\vspace{0.25cm}
This is the case of complex eigenvalues. As $x^*$ is unstable, $\mu >0$. Stabilization is obtained by means of 
\begin{equation} \label{eq12}
   \widetilde{K}=
   \begin{bmatrix} \epsilon_1 &  -\epsilon_2 \\
\epsilon_2 & \epsilon_1  \end{bmatrix}
\end{equation}
where\\
\begin{equation} \label{eq13}
\begin{matrix}
\epsilon_1= \displaystyle\frac{\rho \cos(\theta-\omega\tau)-\rho e^{\mu\tau}\cos(\theta-2\omega\tau)-e^{3\mu\tau}\cos(2\omega\tau)+e^{4\mu\tau}\cos(\omega\tau)}{\tau e^{\mu\tau}(1-2e^{\mu\tau}\cos(\omega\tau)+e^{2\mu\tau})}, \\
\hspace{0.3cm} \\
\epsilon_2=\displaystyle\frac{\rho \sin(\theta-\omega\tau)-\rho e^{\mu\tau}\sin(\theta-2\omega\tau)-e^{3\mu\tau}\sin(2\omega\tau)+e^{4\mu\tau}\sin(\omega\tau)}{\tau e^{\mu\tau}(1-2e^{\mu\tau}\sin(\omega\tau)+e^{2\mu\tau})}, \end{matrix} 
\end{equation} 
\hspace{-0.1cm}for any $\rho$, $\theta$: $0 \leqslant \rho <1$ and $0 \leqslant \theta < 2\pi$. \\

This is obtained by introducing:
\begin{center}
$z=x+iy \ \	\ \wedge \ \ \ 	 \lambda=\mu+i\omega$ \\ 
\end{center}
which yields to the scalar (complex) system:
\begin{equation*}
\dot{z}=\lambda z+ \epsilon (t)\big(z(t-2\tau)-z(t-\tau)\big)
\end{equation*}
where $\epsilon (t)$ results 
as in (\ref{eq7bis}) but with $\epsilon=\epsilon_1+i\epsilon_2$.\\
The methodology of \cite{VEP-GAG} is extrapolated and the map $P$, defined on the complex plane, is associated as in the scalar real case. Then, stabilization is obtained if\\
\begin{equation*}
    e^{3\lambda \tau}[1+\epsilon \tau e^{-2\lambda \tau}(1- e^{\lambda \tau})]= \rho e^{i \theta}
\end{equation*}
and formulae (\ref{eq13}) is deduced.


\section{The n-dimensional case}

Let us assume that the matrix $A$ has $m$ real eigenvalues  $\lambda_1,\dots,\lambda_m$ and $2(n-m)$ complex eigenvalues $\lambda_{m+1},\overline{\lambda_{m+1}},\dots,\lambda_{n},\overline{\lambda_{n}}$, all of them counted  with their multiplicity. Then, the Jordan matrix of $A$ may be written as\\
\begin{equation*} \label{DD}
D = \begin{bmatrix} 
\lambda_1 & 0 		 & \cdots    &  		  &				 &  0 	\\
0 				& \ddots & \ddots	   &        & 			 &			\\
\vdots 		& \ddots & \lambda_m &  			& 			 &  		\\
					&  			 &  				 & B_{m+1}  &  			 &  		\\
					&        & 					 & 				& \ddots &			\\
0					&        & 					 & 				& 			 &	B_n	\\
\end{bmatrix} + 
\begin{bmatrix}  
0 						& \delta_{1,2} 		 & 0    &  \cdots		  &				 &  0 	\\
\vdots 				& \ddots & \ddots	   &        & 			 &			\\
							&  &  &  			& 			 &  		\\
					&  			 &  				 &   &  			 &  		\\
					&        & 					 & 				& \ddots &	\delta_{n-1,n}		\\
0					&        & 					 & 				& 			 &	0	\\
\end{bmatrix} \end{equation*}\\
where $\delta_{ij}=$ 0 or 1 and 
\begin{equation*} B_i = \begin{bmatrix}
\mu_i   &  -\omega_i \\
\omega_i   &  \mu_i \\
\end{bmatrix} 
\end{equation*} 
with $\mu_i$ = Re($\lambda_i$), $\omega_i$ = Im($\lambda_i$) for $i=m+1,\cdots,n$. 

By straightfull generalization of the building of $\widetilde{K}$ in the 1 and 2 dimensional cases, the corresponding matrix $\widetilde{K}$ is determined by\\
\begin{equation*} \widetilde{K}= 
\begin{bmatrix}
k_1  &        &     & 			 &	 \\
     & \ddots &     & 			 &	  \\
     &        & k_m & 			 &	   \\
     &        &     & E_{m+1}&		  \\
     &        &     & 			 &\ddots \\
\end{bmatrix}
\end{equation*}
where\\ $k_i= \begin{cases}
    0 \! & \mbox{if } \lambda_i<0 \\
    \epsilon_i \! & \mbox{if } \lambda_i>0
\end{cases}\ $ and
 $\ E_i= \begin{cases}
    0_{2\times2} \! & \mbox{if  Re}(\lambda_i)<0 \\
    \begin{bmatrix}\epsilon_{i1} & -\epsilon_{i2} \\ \epsilon_{i2} & \epsilon_{i1} \end{bmatrix} \! & \mbox{if  Re}(\lambda_i)>0
\end{cases} \; i=m+1,\dots, n$.\\

\vspace{0.1cm}

\noindent  with $\epsilon_i$ as in (\ref{eq10}) for $i=1,\dots, m$ and $\epsilon_{i1}$ and $\epsilon_{i2}$ as $\epsilon_{1}$ and $\epsilon_{2}$ in (\ref{eq13}) but replacing $\mu, \omega, \rho$ and $\theta$ by $\mu_{i}, \omega_{i}, \rho_{i}$ and $\theta_{i}$ respectively with  
$0\leqslant \rho_{i}<1$, $0 \leqslant \theta_{i}<2 \pi$ for $i=m+1,\dots,n$.\\

Let us note that the computation of $\widetilde{K}$ does not hold on the $\delta_{ij}$ values.\\


\section{ TDFC for controlling chaos}

Let us consider that system (\ref{eq1}) is chaotic. The tools developed in the previous sections may be implemented for controlling chaos by introducing an strategy based on the algorithm (\ref{eq3})-(\ref{eq3bis}). Let $\delta$ be such that for all initial condition within a $\delta$-distance from the equilibrium $x^*$, 
asymptotic stabilization in $x^*$ is achieved by applying (\ref{eq3})-(\ref{eq3bis}). The strategy consists in two steps. In the first stage, the control is not activated (i.e., the system runs free). The second stage begins at a certain time $t^*$ for which the trajectory $x(t)$ satisfies $||x(t)-x^*||< \delta$ and control is activated. The existence of $t^*$ (usually called ``wait time") is due to the ergodicity of the chaotic system (\cite{BOCCALETTI2000103}) and it is computed on line. Namely 

\begin{equation} \label{eqControl}
   u_{c}(t)= \begin{cases}  0 \! & \mbox{if }0 \leqslant t \leqslant t^* \\  
	 K_c(t)\big(x(t-2\tau)-x(t-\tau)\big) & \mbox{if }  t \geqslant t^* \end{cases}
\end{equation}
where $t^*=inf\{t\geqslant 2\tau, ||x(t)-x^*||\leqslant\delta\}$ and $K_c(t)=K(t-t^*)$ for $K(t)$ defined in (\ref{eq3bis}).\\
This strategy applies if the equilibrium point to stabilize is embedded in the attractor. Its effective implementation requires some knowledge on the basin of attraction of the attractor, i.e., a trapping region of it. \\
For illustration, the strategy (\ref{eqControl}) will be used to control chaos of the classical chaotic systems known as Chua's and Rossler systems in some of its equlibrium.


\subsection{Chua's equation}

Chua's circuit is a classical device to study chaos phenomena. The circuit resulting from adding a linear resistor in series with the inductor has been named as Chua's oscillator which presents different dynamical features for certain variations on its parameters (\cite {Galeria}). 

Its dimensionless equation is given by: 
\begin{equation} \label{eqIVdeG}
\begin{cases}
\dot{x_1}=\alpha\big(x_2-x_1-f(x_1)\big)\\
\dot{x_2}=x_1-x_2+x_3 \\
\dot{x_3}=-\beta x_2-\gamma x_3
\end{cases}
\end{equation} 
where,
\begin{equation} \label{eqVdeG}
f(x_1)=m_1x_1+\frac{1}{2}(m_0-m_1)\big(|x_1+1|-|x_1-1|\big)
\end{equation} 

Let us note that (\ref{eqIVdeG})-(\ref{eqVdeG}) may be re-written as the following piece-wise linear differential equation
\begin{equation} \label{eq39-rgsequi}
   \dot{x}= \begin{cases}  A_1 x(t)-d \! & \mbox{if } x \in D_{-1}= \left\{(x,y,z):x \leqslant -1\right\} \\
	A_2 x(t) \! & \mbox{if } x \in D_0 = \left\{(x,y,z):\left|x\right| \leqslant 1\right\} \\
	A_1 x(t)+d\! & \mbox{if } x \in D_{1}= \left\{(x,y,z):x \geqslant 1\right\} \\\end{cases}
\end{equation}
being 
\begin{equation*}
A_1=\begin{bmatrix} -\alpha(1+m_1) & \alpha & 0\\ 1 & -1 & 1 \\ 0 & -\beta & -\gamma \end{bmatrix}, \  
A_2=\begin{bmatrix} -\alpha(1+m_0) & \alpha & 0\\ 1 & -1 & 1 \\ 0 & -\beta & -\gamma \end{bmatrix},\  
d=\begin{pmatrix}\alpha (m_0-m_1) \\0 \\0\end{pmatrix}
\end{equation*} 
 
In the general case, system (\ref{eq39-rgsequi})
presents three equilibrium points and each one belongs to each $D_i$. In particular, the origin is one of them.  By varying the parameters set, Chua's equation presents different analytical and geometric features (\cite{Galeria}).

Chua's oscillator reduces to Chua's circuit upon setting the linear resistor $R_0$ to zero, i.e., $\gamma =0$ which includes the cases known as double scroll and the double hook because of the aspect of their respective attractors. For both of them, the origin is in the attractor and it is an equilibrium point having only one positive eigenvalue, although the linearization in the origin of the double scroll has one real eigenvalue and two complex eigenvalues while for the double hook has three real eigenvalues.

In \cite{GALIAS2022106109} rigourous ad-hoc algorithms are developed to compute solutions of continuous piece-wise linear systems. These tools are used to obtain trapping regions enclosing the scroll and hook attractors which are required to implement the control strategy (\ref{eqControl}) on the mentioned Chua's circuits. Besides, these tools are effective to compute trajectories passing too close to the origin, a feature to take into account if (\ref{eqControl}) is to be applied with $\delta$ very small.


\subsubsection{Double scroll stabilized in the origin}
The study on the double scroll attractor dates back to the middle eighties (\cite{Matsu85}) and it was sustained by a great deal of work afterwards but a rigorous analysis of its complex structure remains open (\cite{Bely2023}). The case  given by (\ref{eq39-rgsequi}) with $\alpha= 9.3515$, $\beta=14.79$, $\gamma=0$, $m_0=-1.138$ and $m_1=-0.722$ is considered in \cite{GALIAS2022106109} where  a trapping region given by the union of four polygons is identified and trajectories starting in the border of one of these polygons are computed. In particular, it is shown how close to the origin the trajectory initiated in $CI=(1.5,\ -0.253849008275,\ -2.55651050226)$ passes. 

The objective is to control chaos by stabilizing the system in the origin. The linearization matrix is given by $A_2$ which has the eigenvalues: $\lambda_1 \approx 2.2407$, $\lambda_{2,3} \approx -0.9751 \pm 2.7508i$. Hence the origin is a saddle-focus equilibrium.

The control strategy (\ref{eqControl}) will be implemented with  $\zeta=0.4$. Taking  $\tau=0.1$ the following matrix $K$ is obtained:

\begin{equation*}
K=\begin{bmatrix} 46.8954 & 123.5023 & 14.6824 \\ 4.7652 & 12.5495 & 1.4919 \\ -31.4526 & -82.8327 & -9.8475 \end{bmatrix}
\end{equation*}

Putting $\delta=1.8$, the results from simulation of the controlled system with $CI$ as initial condition are appreciated in 
Figures \ref{Fig1-DS-A} and \ref{Fig2-DS-EC}.  The free and controlled system initiated in $CI$ are depicted in Figure \ref{Fig1-DS-A}(a) and Figure \ref{Fig1-DS-A}(b) respectively, for confontation. The $\delta$-sphere centered in the origin depicted in both graphs put in evidence the effect of the control on the states bound. It may be seen in Figure \ref{Fig1-DS-A}(b)  that once the control is activated, the trajectory does not leave the sphere of ratio $\delta$  any more. State behaviour and control performance is also appreciated by the state and control signals of the controlled system  displayed in Figure \ref{Fig2-DS-EC}.

\begin{figure}[h!]
    \centering
    \includegraphics[width=0.49\textwidth]{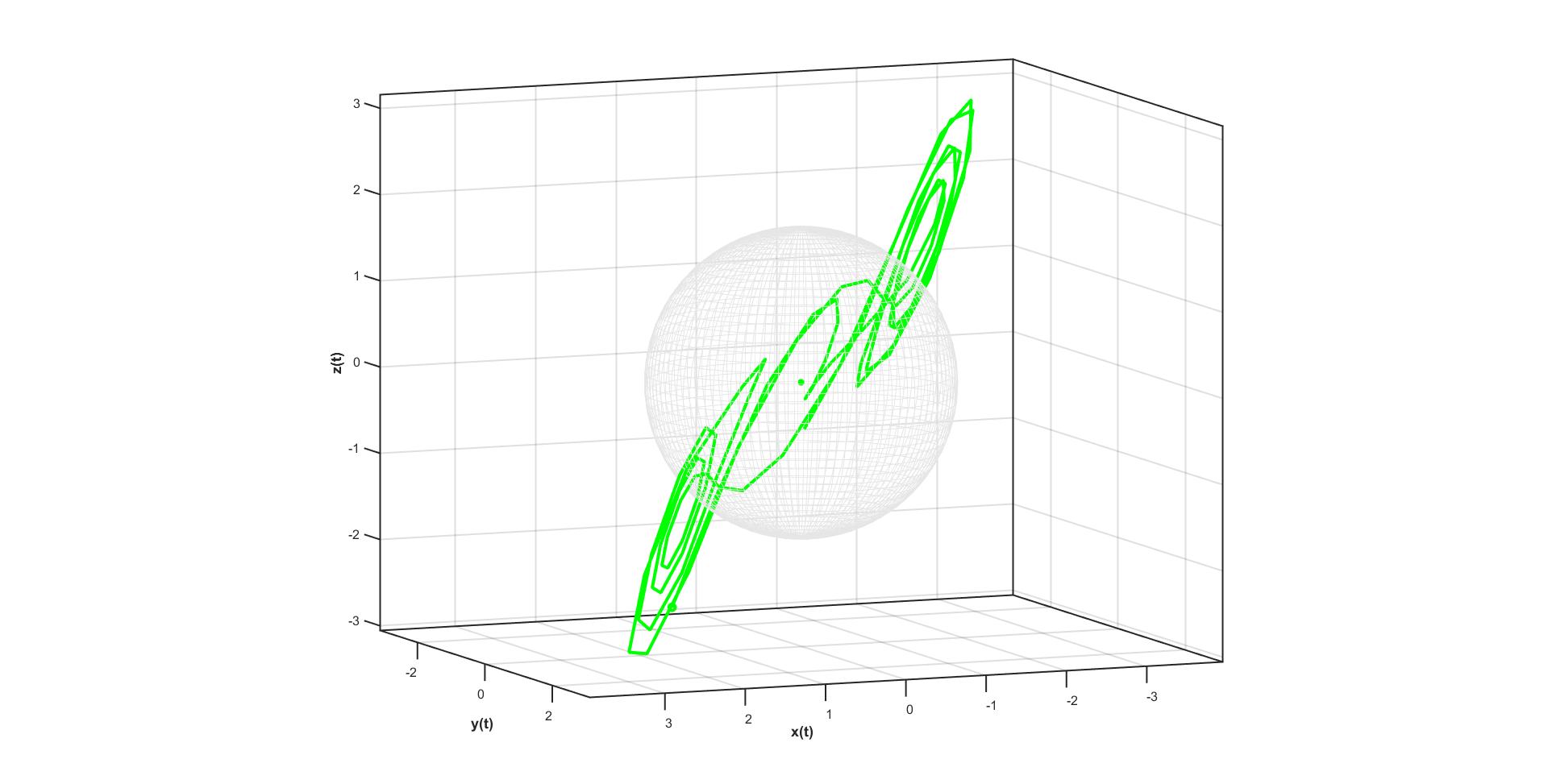}    
    \includegraphics[width=0.49\textwidth]{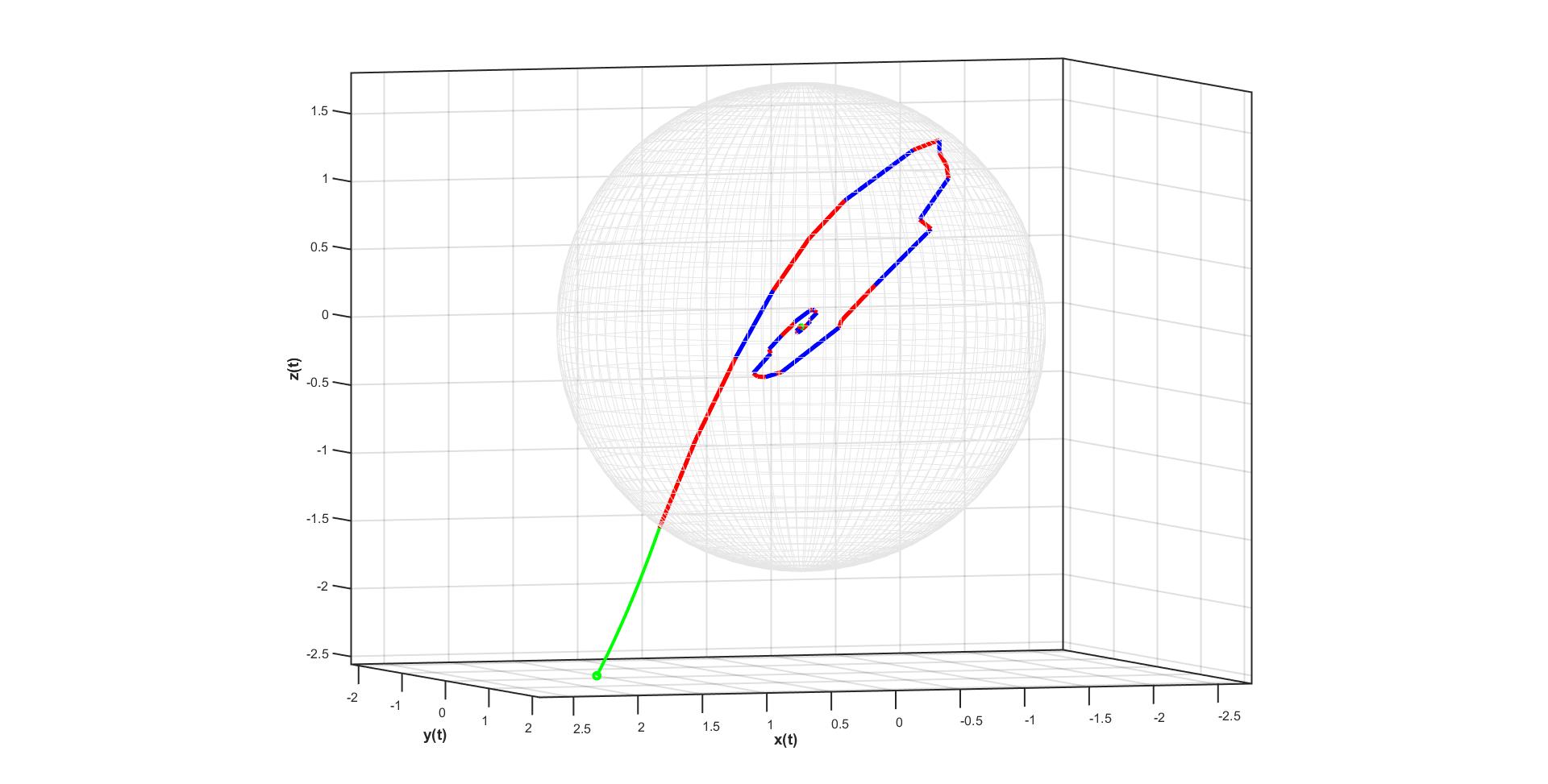}
    \caption{Double scroll. Trajectory initialized in $CI=(1.5,-0.253849008275,-2.55651050226)$ of a) the free and b) the controlled system with $\zeta = 0.4$, $\tau = 0.1$ and $\delta=1.8$.}
    \label{Fig1-DS-A}
\end{figure}
\begin{figure*}[h!]
    \centering
    \includegraphics[width=0.5\linewidth, height=0.15\textheight]{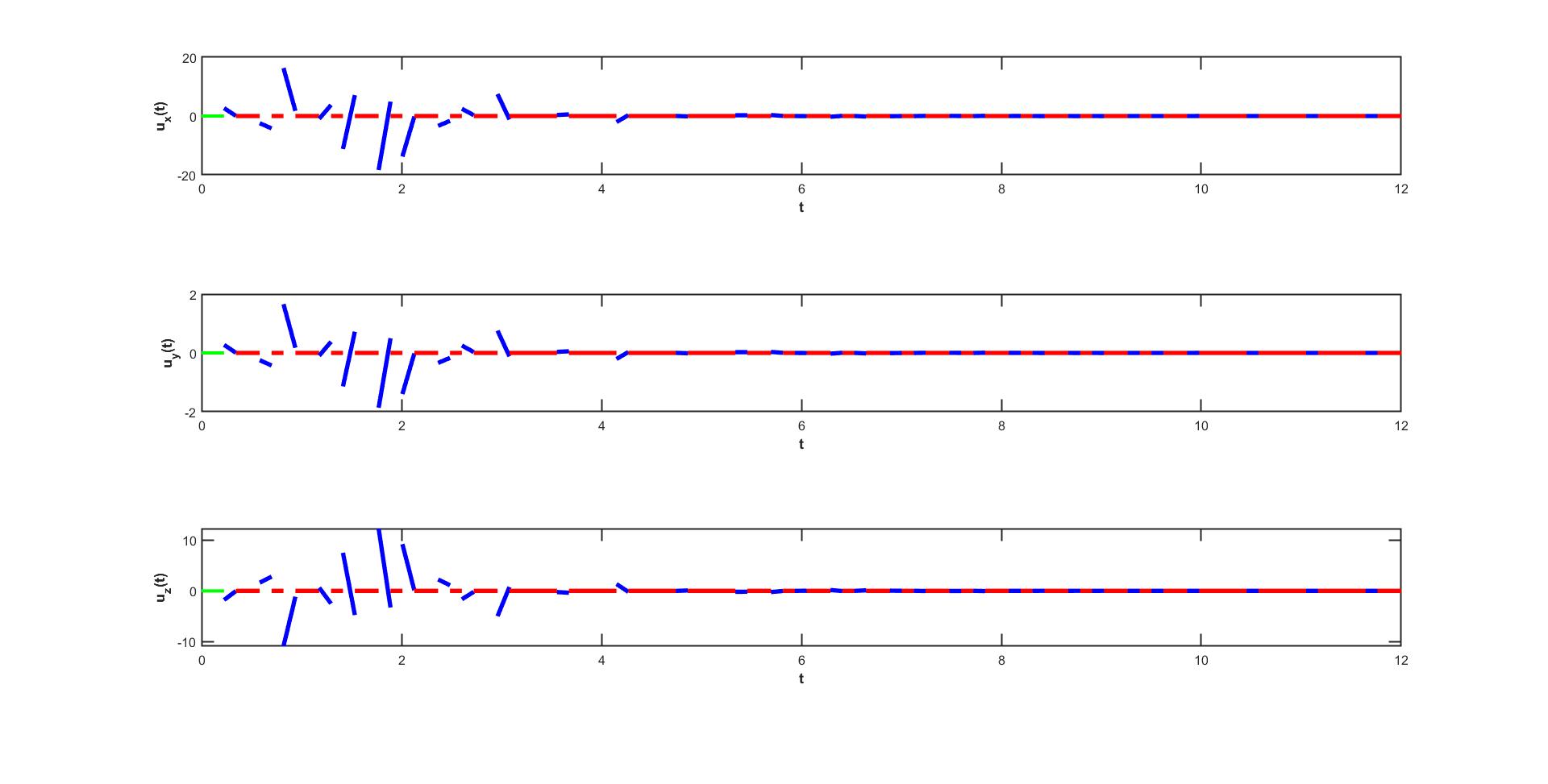}%
    \includegraphics[width=0.5\linewidth, height=0.15\textheight]{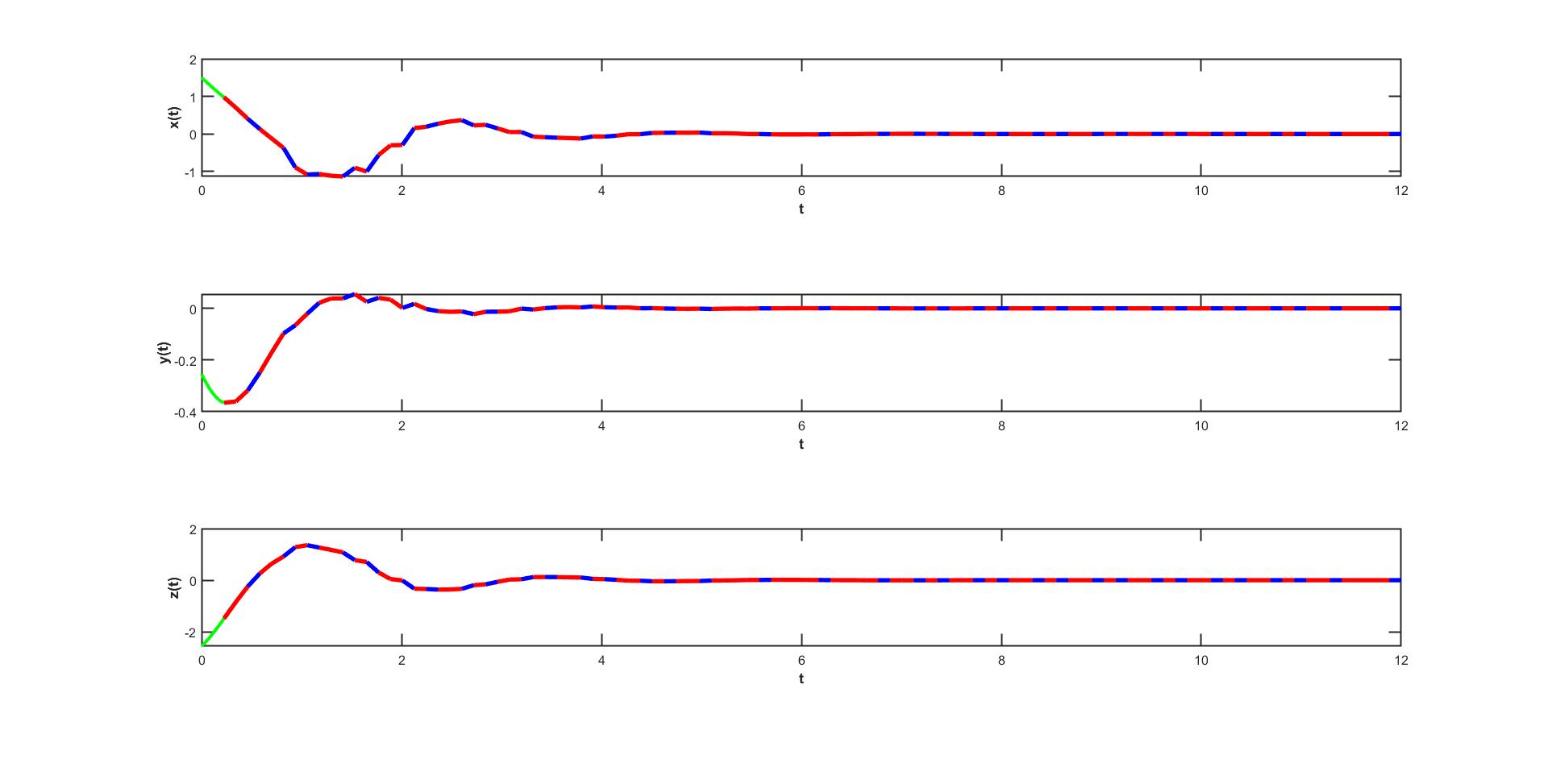}
    \caption{Double scroll. State  and control signals of the controlled system initialized in $CI=(1.5,-0.253849008275,-2.55651050226)$  with $\zeta = 0.4$, $\tau = 0.1$ and $\delta=1.8$.}
    \label{Fig2-DS-EC}
\end{figure*}

Simulation of the controlled system, keeping the $\zeta$ and $\tau$ values but choosing $\delta=0.5$ is illustrated in Figure \ref{Figura 3} where the resulting trajectory in the three dimensional space and the respective temporal states and control signals are shown.  Let us note that although the trajectory goes out the sphere during some time, it eventually remains there and achieves convergence.

\begin{figure}[h!] 
    \centering
    \begin{minipage}{.5\textwidth}
        \centering
        \includegraphics[width=1\linewidth]{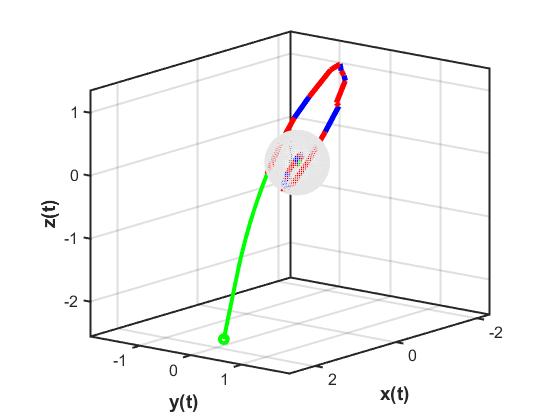}
    \end{minipage}%
    \begin{minipage}{0.5\textwidth}
        \centering
        \includegraphics[width=1\linewidth, height=0.15\textheight]{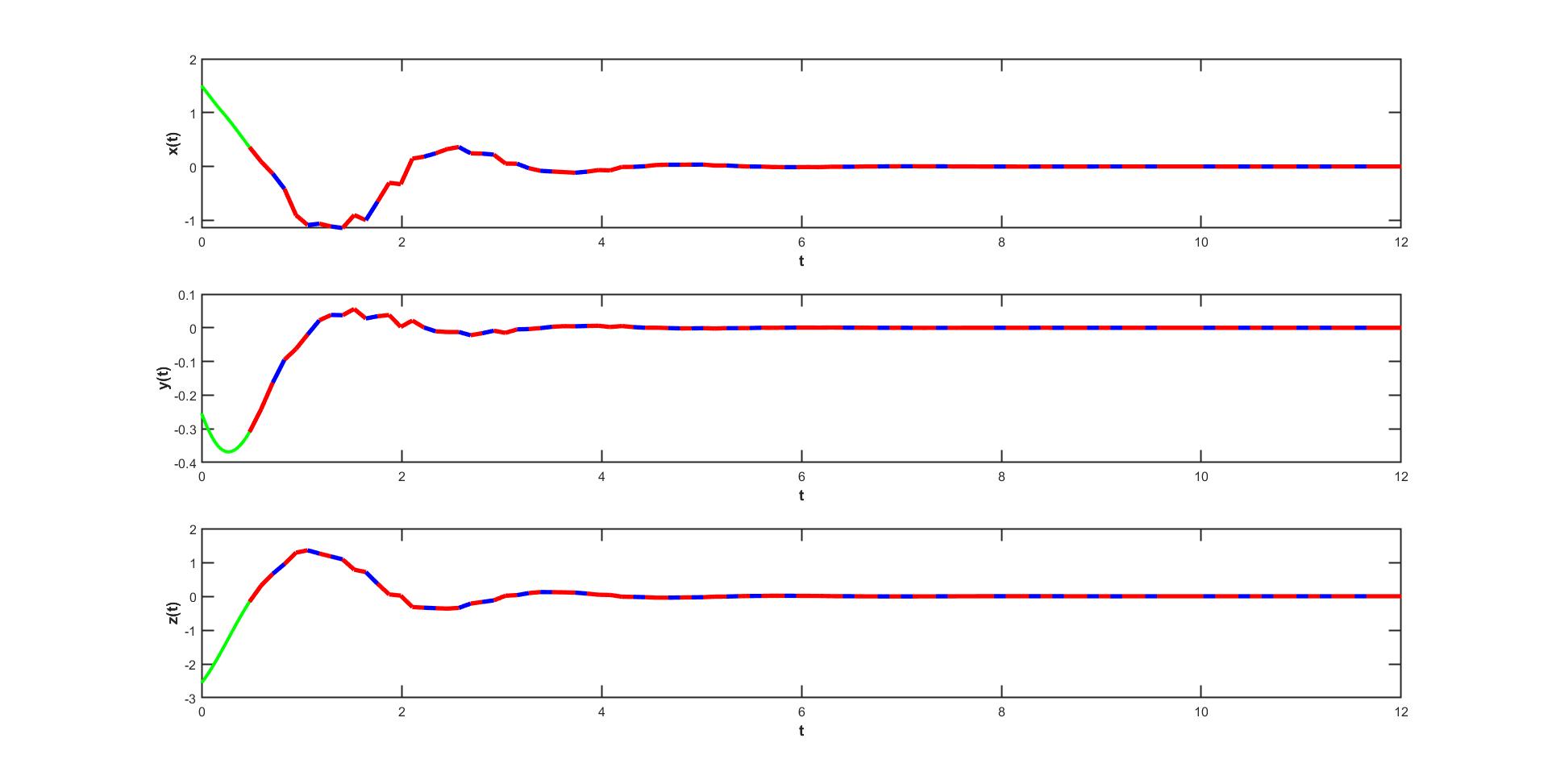}\\
        \includegraphics[width=1\linewidth, height=0.15\textheight]{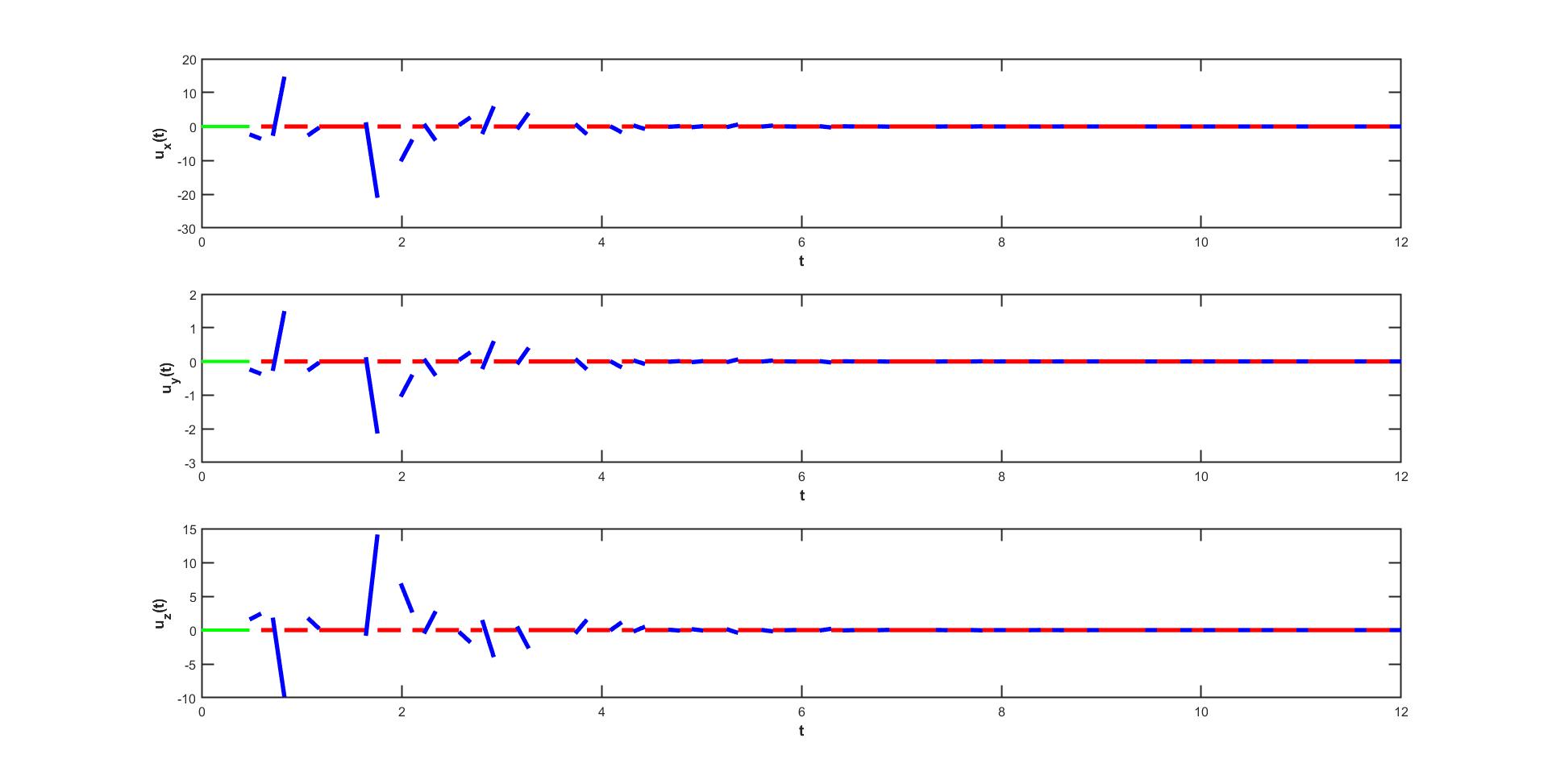}        
    \end{minipage}
    \caption{Double scroll controlled with $\zeta = 0.4$, $\tau = 0.1$ and $\delta=0.5$. Trajectory initiated in  $CI=(1.5,\ -0.253849008275,\ -2.55651050226)$: a) in 3D-space, b) state and control signals.}
    \label{Figura 3}
\end{figure}


\subsubsection{Double hook stabilized in the origin.}

The double hook attractor was presented in \cite{9914} when introducing a totally different parameter set and it goes on being a subject of interest (\cite{Galias2022}). The case considered in  \cite{GALIAS2022106109} belongs to (\ref{eq39-rgsequi}) with $\alpha= -6$, $\beta=-4.442$, $\gamma=0$, $m_0=-2.265$ and $m_1=-0.93$.  As in the previous example, a trapping region  is computed in \cite{GALIAS2022106109} and  in particular, the point $CI=(18,\ 1.5776491929,\ -14.3838522424)$ belongs to it.

The objective is to control chaos by stabilizing the system in the origin. The linearization matrix is given by $A_2$ which has the eigenvalues: $\lambda_1 \approx 1.4336,\ \lambda_2 \approx -3.7467,\ \lambda_3 \approx -6.2768$. Hence, the origin is a saddle equilibrium.

The control is designed by taking $\zeta=0.4$ and $\tau=0.25$ yielding to the following matrix $K$:

\begin{equation*}
K=\begin{bmatrix}
   11.4871 & 16.4678 & 1.8250 \\
   3.7073 & 5.3147 & 0.5890 \\
   -2.4651 & -3.5339 & -0.3916
\end{bmatrix}
\end{equation*}  

The control strategy (\ref{eqControl}) is implemented on this system with $\delta=6$ and it is simulated with initial condition $CI$. Confrontation of the free and the controlled systems are appreciated in Figure \ref{Fig4-DH-A}, showing the respective trajectories in relation to the sphere of ratio 6. In Figure \ref{Fig4-DH-A}(b), it can be seen that the trajectory evolves freely enough time  to reveal the double hook structure. State behaviour and control performance are also appreciated through the evolution of each signal in Figure \ref{Fig4-DH-EC}.
Convergence is achieved even activating the control not so close from the origin, as it happens for the same case but with $\delta=10$ (see Figure \ref{Figure 6}) although  affecting the control performance, i.e. the time of convergence is reduced but in turn, it demands more control effort.

\vspace{1cm}

\begin{figure*}[h!]
  \includegraphics[width=0.45\textwidth]{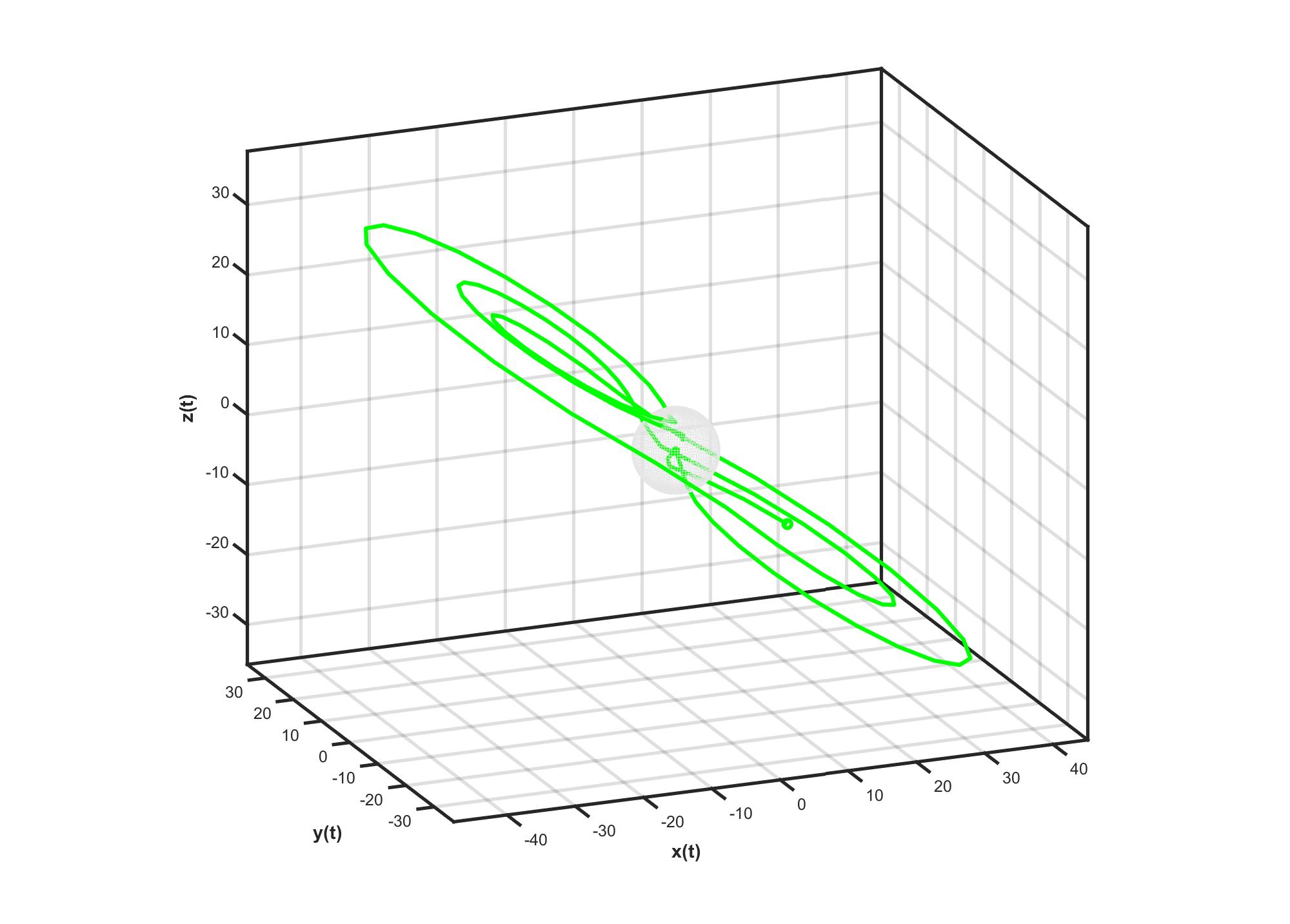} 
  \includegraphics[width=0.65\textwidth]{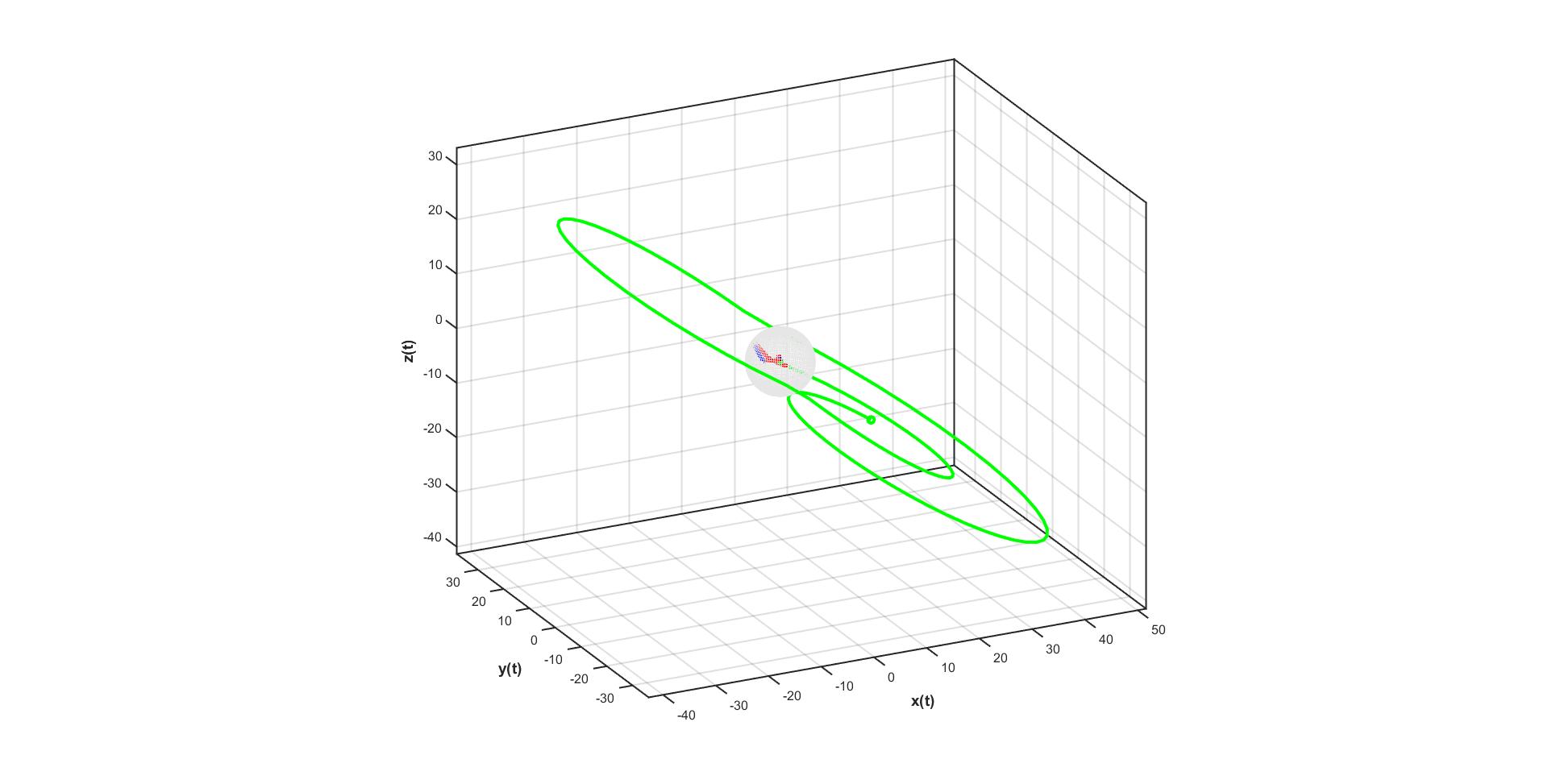}
\caption{Double hook. Trajectory initialized in $CI=(18,\ 1.5776491929,\ -14.3838522424)$ of a) the free and b) the controlled system with $\zeta = 0.4$, $\tau = 0.25$ and $\delta=6$.}  
\label{Fig4-DH-A} 
\end{figure*}

\begin{figure*}[h!]  
  \includegraphics[width=0.5\linewidth, height=0.15\textheight]{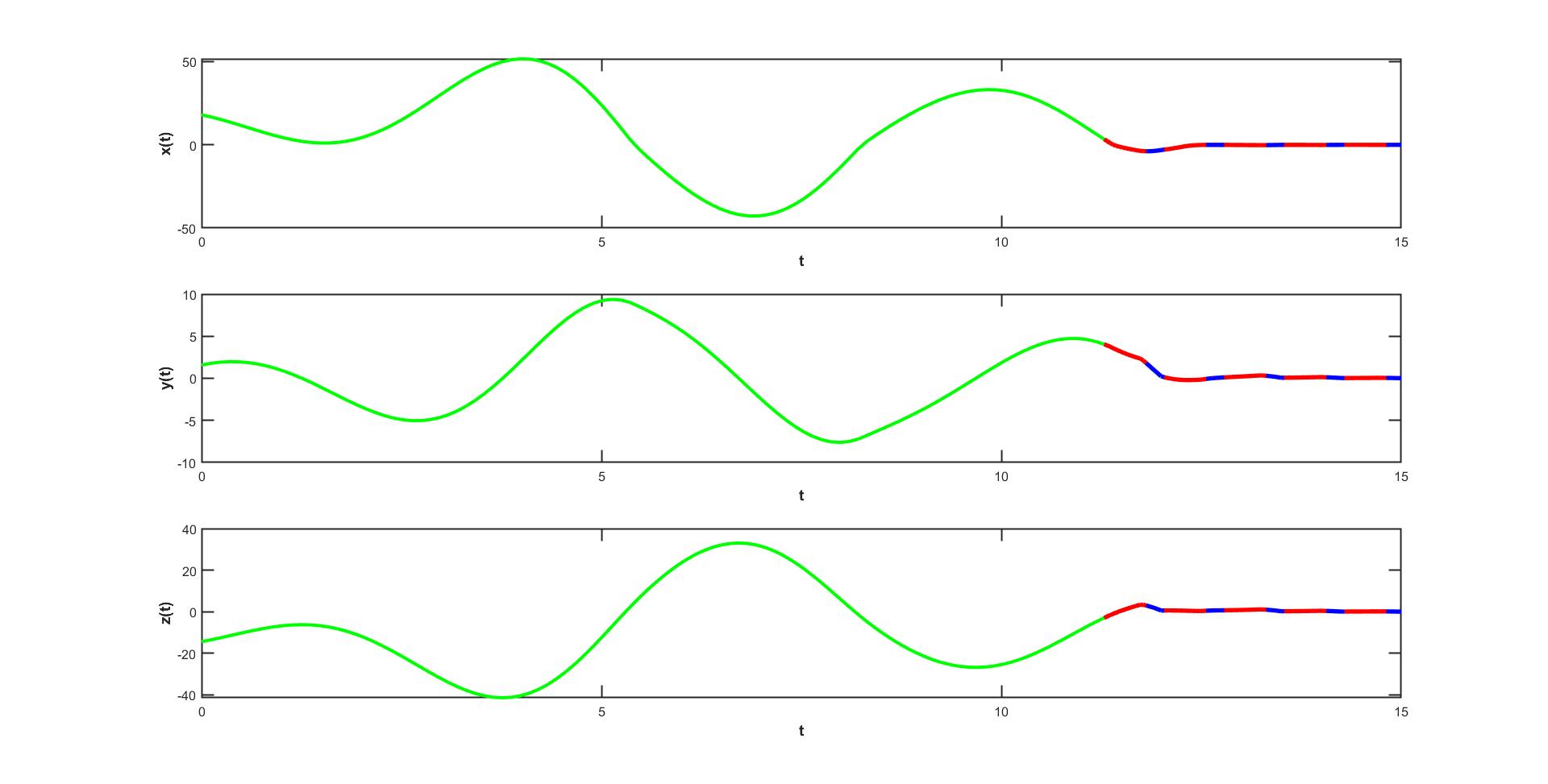}
  \includegraphics[width=0.5\linewidth, height=0.15\textheight]{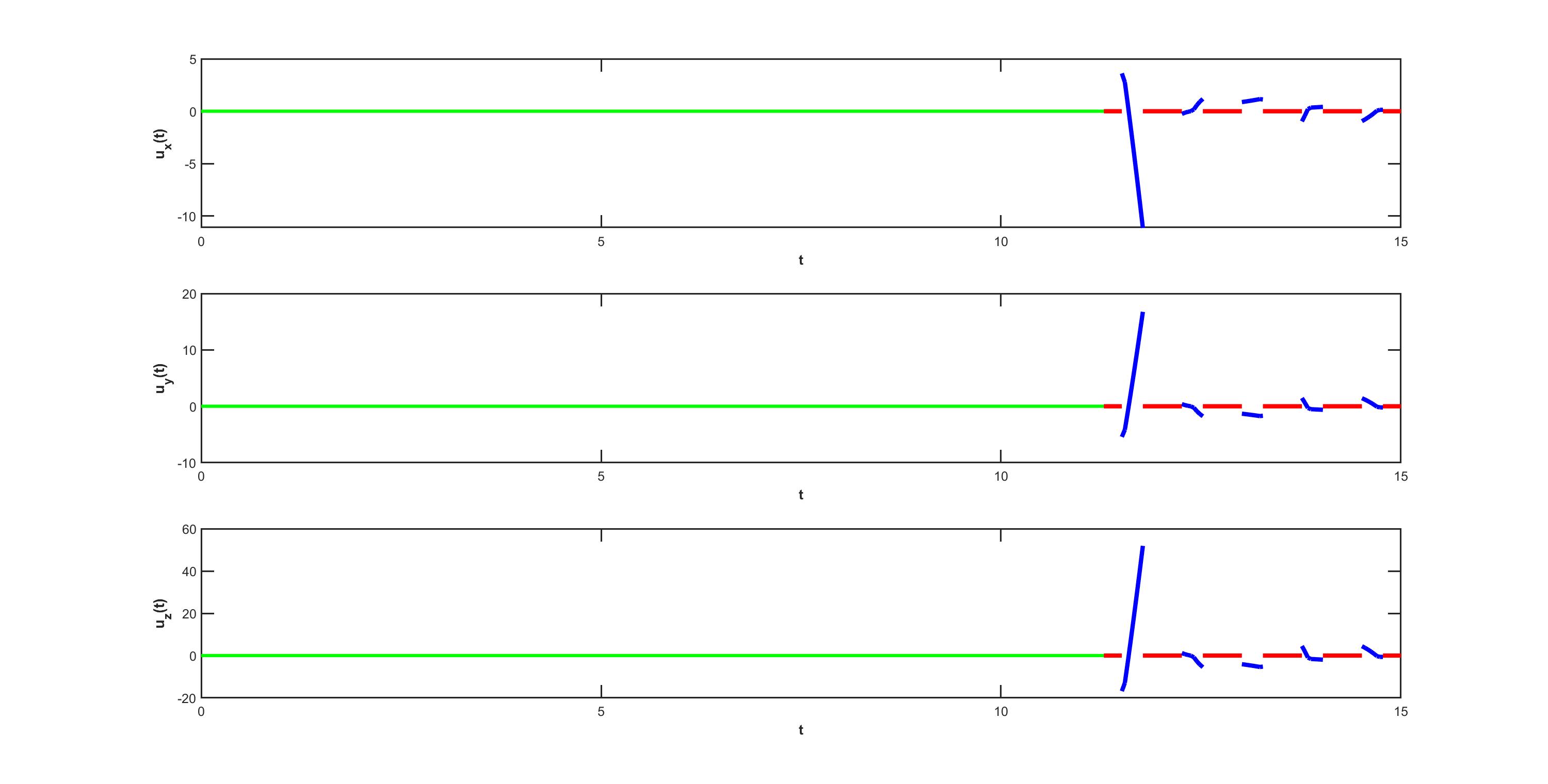}
\caption{Double hook. State  and control signals of the controlled system initialized in $CI=(18,\ 1.5776491929,\ -14.3838522424)$  with $\zeta = 0.4$, $\tau = 0.25$ and $\delta=6.$}
\label{Fig4-DH-EC}       
\end{figure*}

\begin{figure}[!htb]
    \centering
    \begin{minipage}{.5\textwidth}
        \centering
        \includegraphics[width=1.2\linewidth]{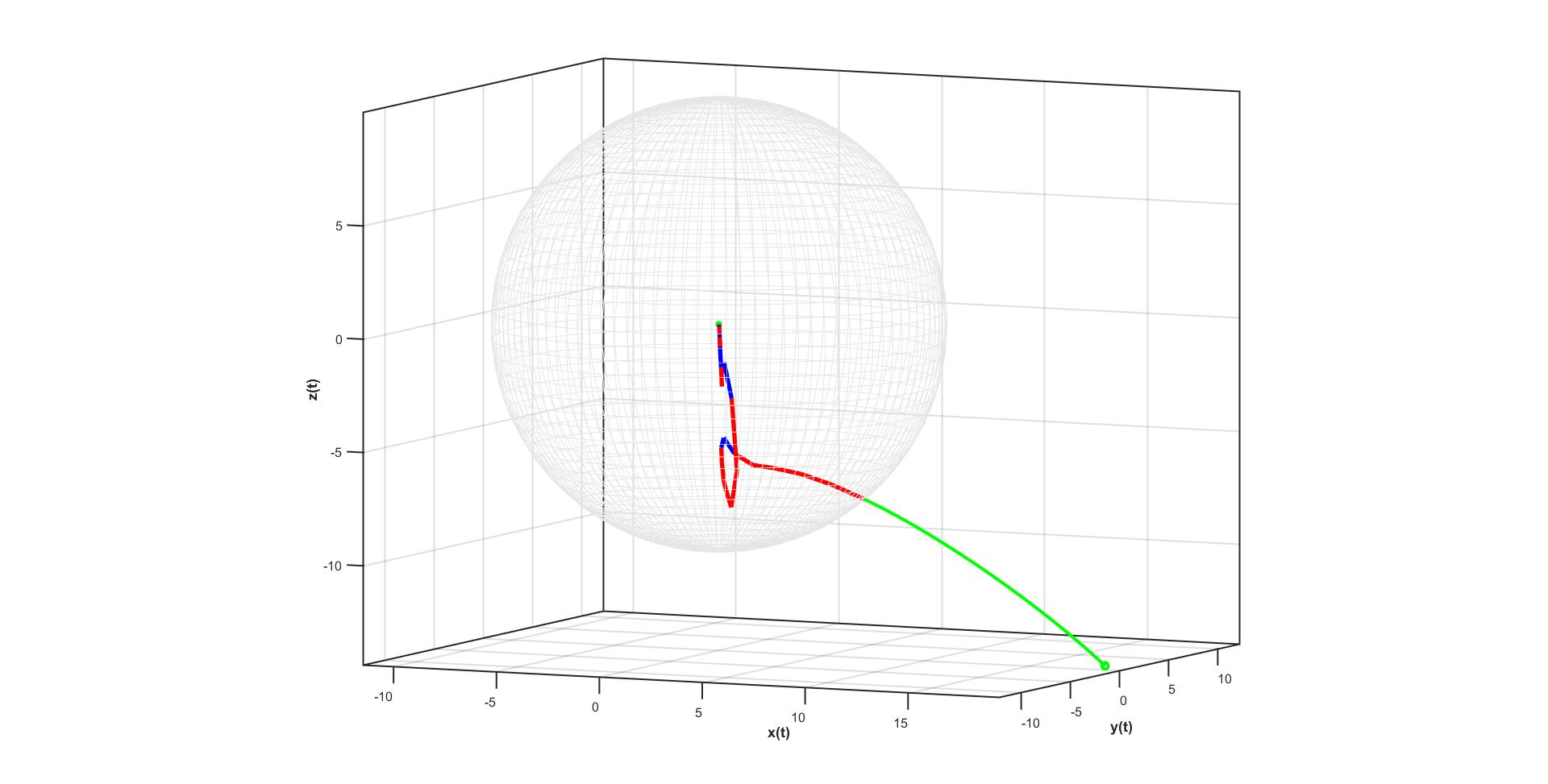}
    \end{minipage}%
    \begin{minipage}{0.5\textwidth}
        \centering
        \includegraphics[width=1\linewidth, height=0.15\textheight]{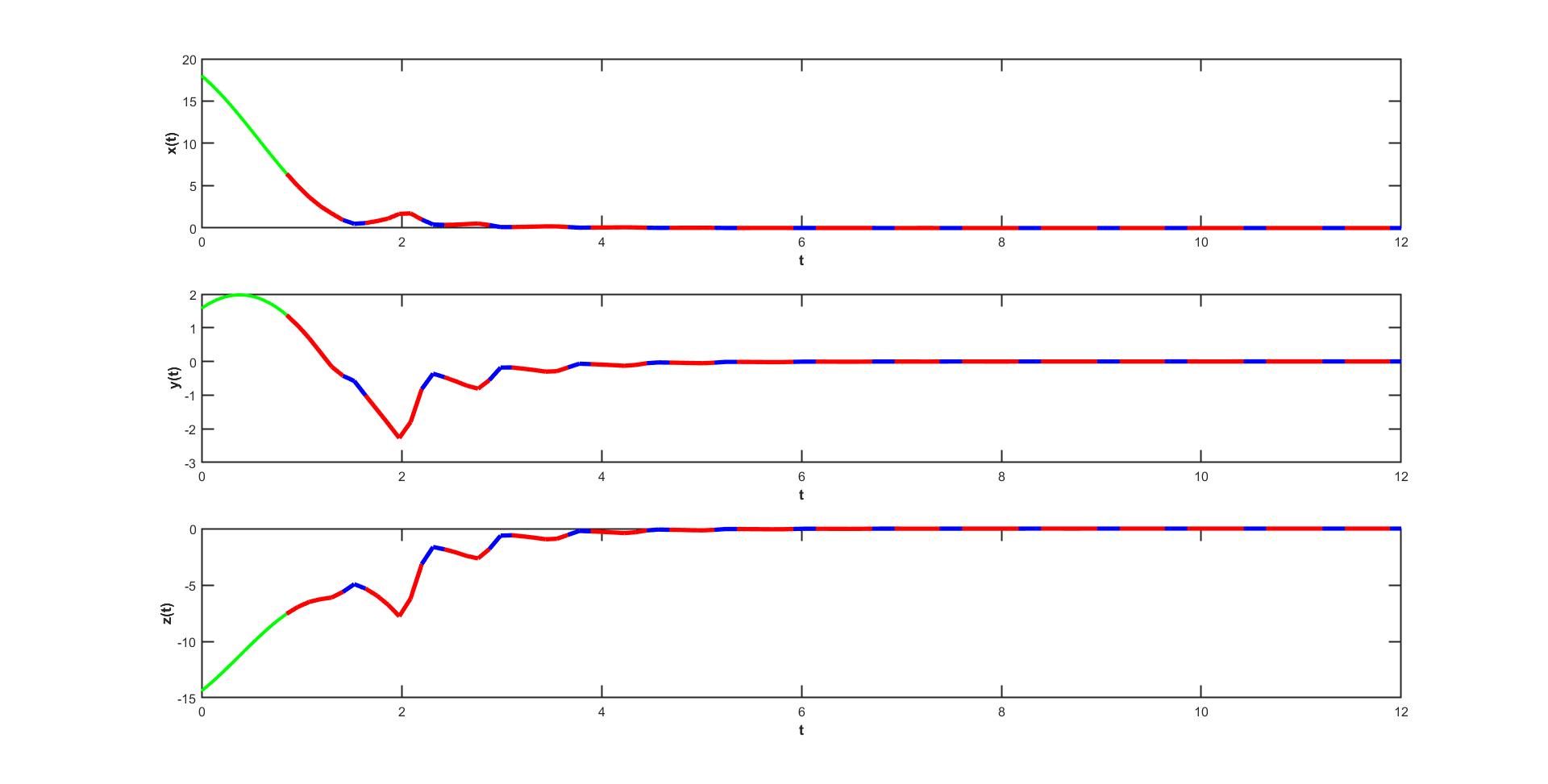}\\
        \includegraphics[width=1\linewidth, height=0.15\textheight]{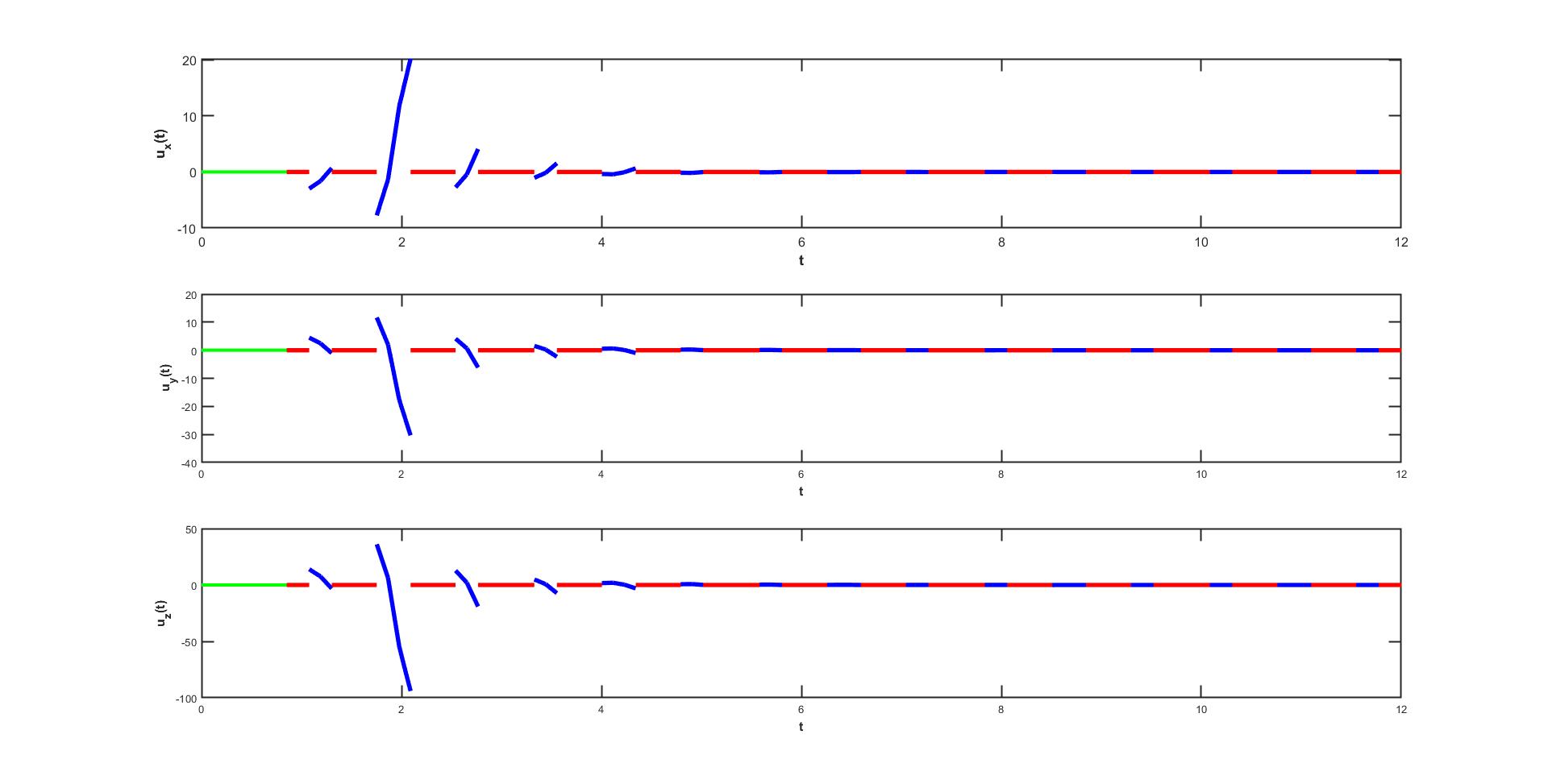}        
    \end{minipage}
    \caption{Double hook controlled with $\zeta = 0.4$, $\tau = 0.25$ and $\delta=10$. Trajectory initiated in  $CI= (18,\ 1.5776491929,\ -14.3838522424)$: a) in 3D-space, b) state and control signals.}
    \label{Figure 6}
\end{figure}

\textit{Remark}: In the both previous examples the linearization matrix has only one positive eigenvalue so the technique of \cite{Lenov2018} also applies but the proposal of the present work involves less parameters and besides the set of admissible values is fully described.

\subsubsection{Example 3: Double scroll stabilized out of the origin.}

The double scroll case is revisited but for stabilizing in another of its equilibrium points. This example is taken from \cite{Konishi2011} and it coincidences with (\ref{eq39-rgsequi}) for $\alpha=9,\ \beta=100/7,\ \gamma=0,\ m_0=-8/7$ and  $m_1=-5/7$. In \cite{Konishi2011}, the control objective is achieved by applying a delayed feedback control based on the act-and-wait concept.

The objective is to control chaos by stabilizing the system in the equilibrium point $P=(-1.5, 0, 1.5)$ which linearization matrix   has the following three eigenvalues: $\lambda_1 \approx -3.9421,\ \lambda_{2,3} \approx 0.1854 \pm 3.0470i$.

Then, by choosing $\rho=0.6$, $\theta=\arctan (\omega/\mu) \approx 1.5100$ where $\mu=$ Re$(\lambda_2)$; $\omega=$Im$(\lambda_2)$ and, $\tau=0.25$, the matrix $K$  results: 
\begin{equation}\label{K-DSfq}
K=\begin{bmatrix} -0.9088 & 0.0792 & -1.6685 \\ 1.0841 & 5.1815 & 0.5345 \\ 1.3856 & -5.3236 & 3.9796\end{bmatrix}
\end{equation}  

The control strategy (\ref{eqControl}) is applied  by taking $\delta=3$ and stabilization is achieved  as it can be seen in Figure \ref{Fig7-DSFE}. The trajectory of the controlled system  with initial condition $CI=(2.104,\ -0.3188,\ -2.0866)$ in phase space is displayed in \ref{Fig7-DSFE}(a) while the state and control signals are appreciated in \ref{Fig7-DSFE}(b).

\begin{figure}[!htb]
    \centering
    \begin{minipage}{.5\textwidth}
        \centering
        \includegraphics[width=1.3\linewidth]{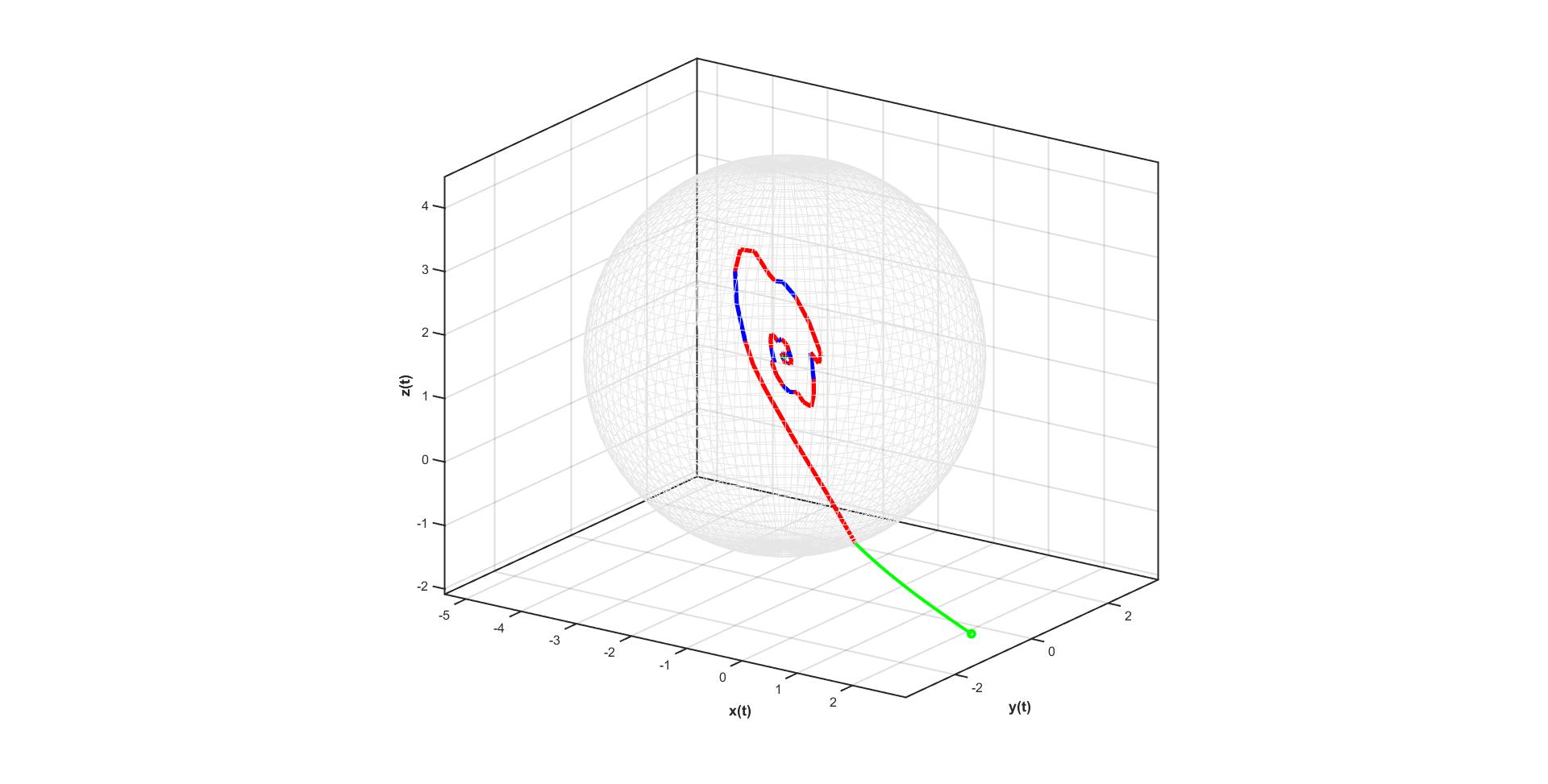}
    \end{minipage}%
    \begin{minipage}{0.5\textwidth}
        \centering
        \includegraphics[width=1\linewidth, height=0.15\textheight]{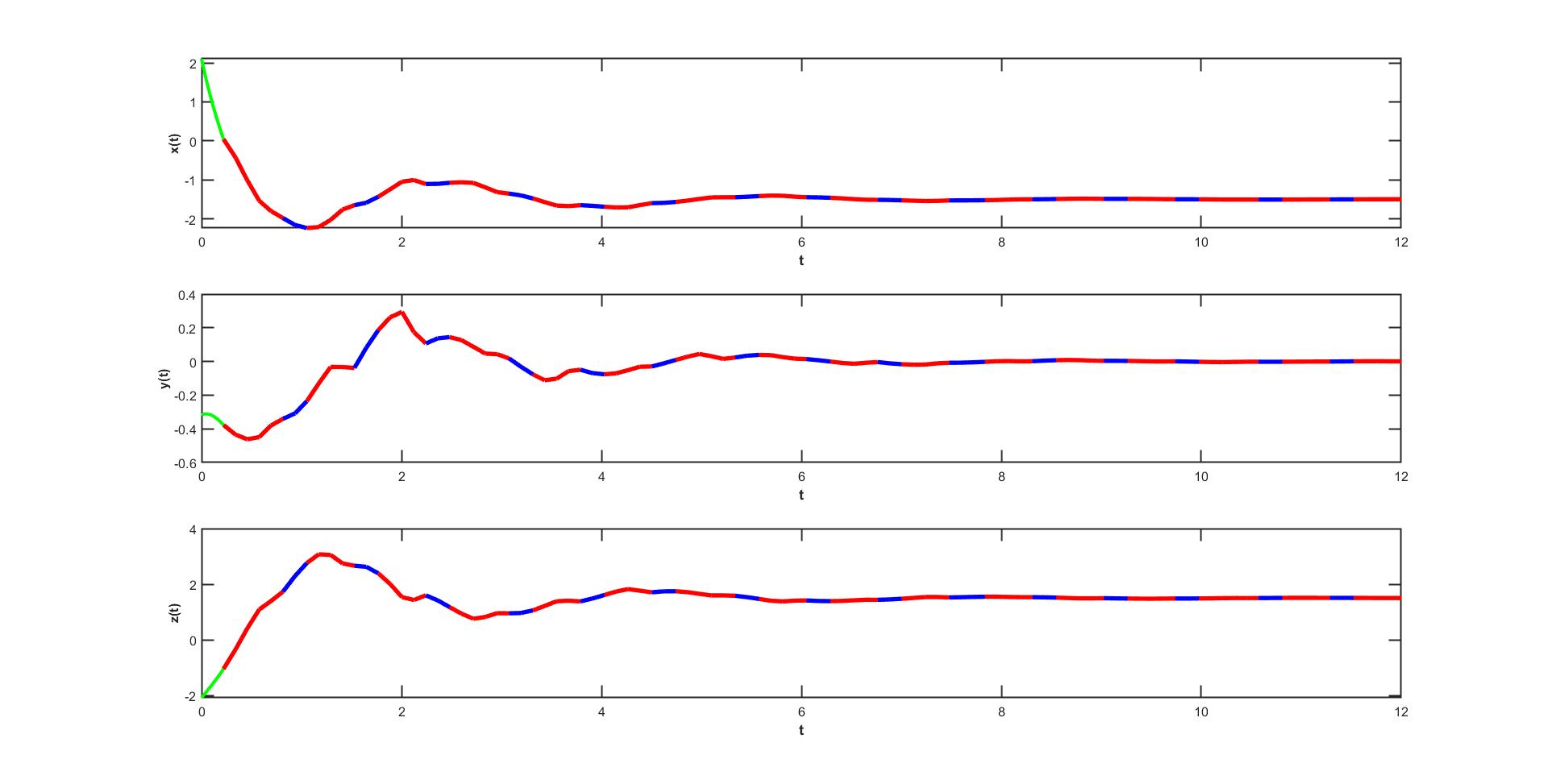}\\
        \includegraphics[width=1\linewidth, height=0.15\textheight]{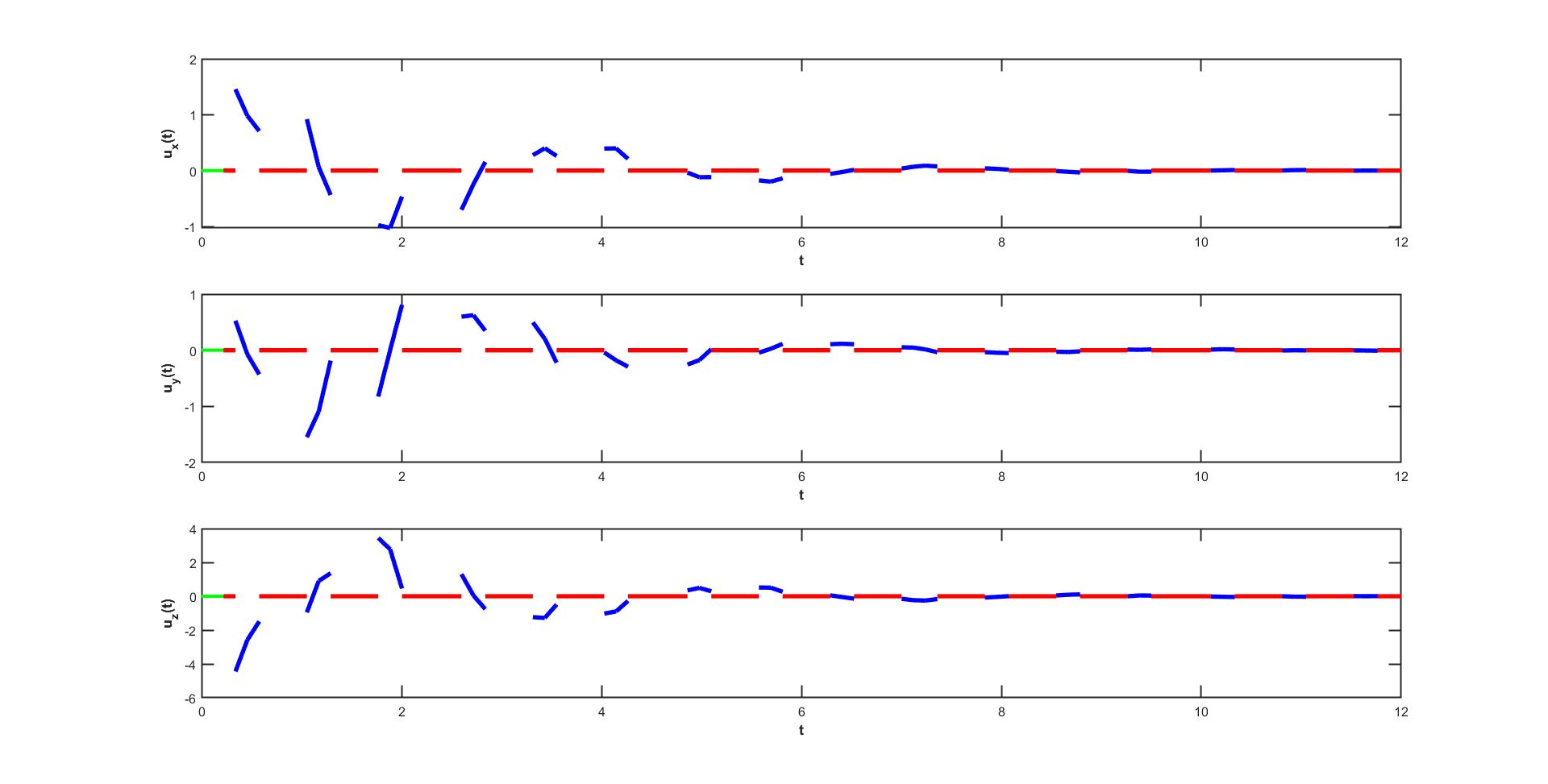}        
    \end{minipage}
    \caption{Double scroll controlled with $\rho=0.6$, $\theta \approx 1.5100$, $\tau = 0.25$ and $\delta=3$. Trajectory initiated in  $CI=(2.104,\ -0.3188,\ -2.0866)$: a) in 3D-space, b) state and control signals.}
    \label{Fig7-DSFE}
\end{figure}

\medskip

Remark: This case is out of the ONL and hence, the methodology based on the act-and-wait concept also applies (\cite{Konishi2011}). The bounds of the state and control signals obtained by using this method or the strategy (\ref{eqControl}) are similar (confront Figure 5 of \cite{Konishi2011} and \ref{Fig7-DSFE}(b)). As the initial condition is not near enough to the equilibrium point, for the first method the state signals (see Figures 5 and 6 of \cite{Konishi2011}) leaves oscillatory behaviour too much later than by means of (\ref{eqControl}) (see Figure \ref{Fig7-DSFE}(b)); this is thank to the introduction of the parameter $\delta$ (and choosing an adequate value for it).

\medskip

It is worth to point out the importance of choosing the appropriate parameter values. Indeed, the control objective is not achieved only by introducing a little change in the delay value; namely fixing $\tau=0.29$, the matrix $K$ comes out:
\begin{equation}\label{K-DSfq}
K=\begin{bmatrix}
0.6831 &	5.9792 &	-0.4122 \\
1.6423 &	6.7941 &	1.1008 \\
-3.4647 &	-22.8502 &	0.0278 \end{bmatrix}
\end{equation}  
but the resulting trajectory diverges as shown in Figure \ref{Figure 8}.

\begin{figure}[!htb]
    \centering
    \begin{minipage}{.5\textwidth}
        \centering
        \includegraphics[width=1.3\linewidth]{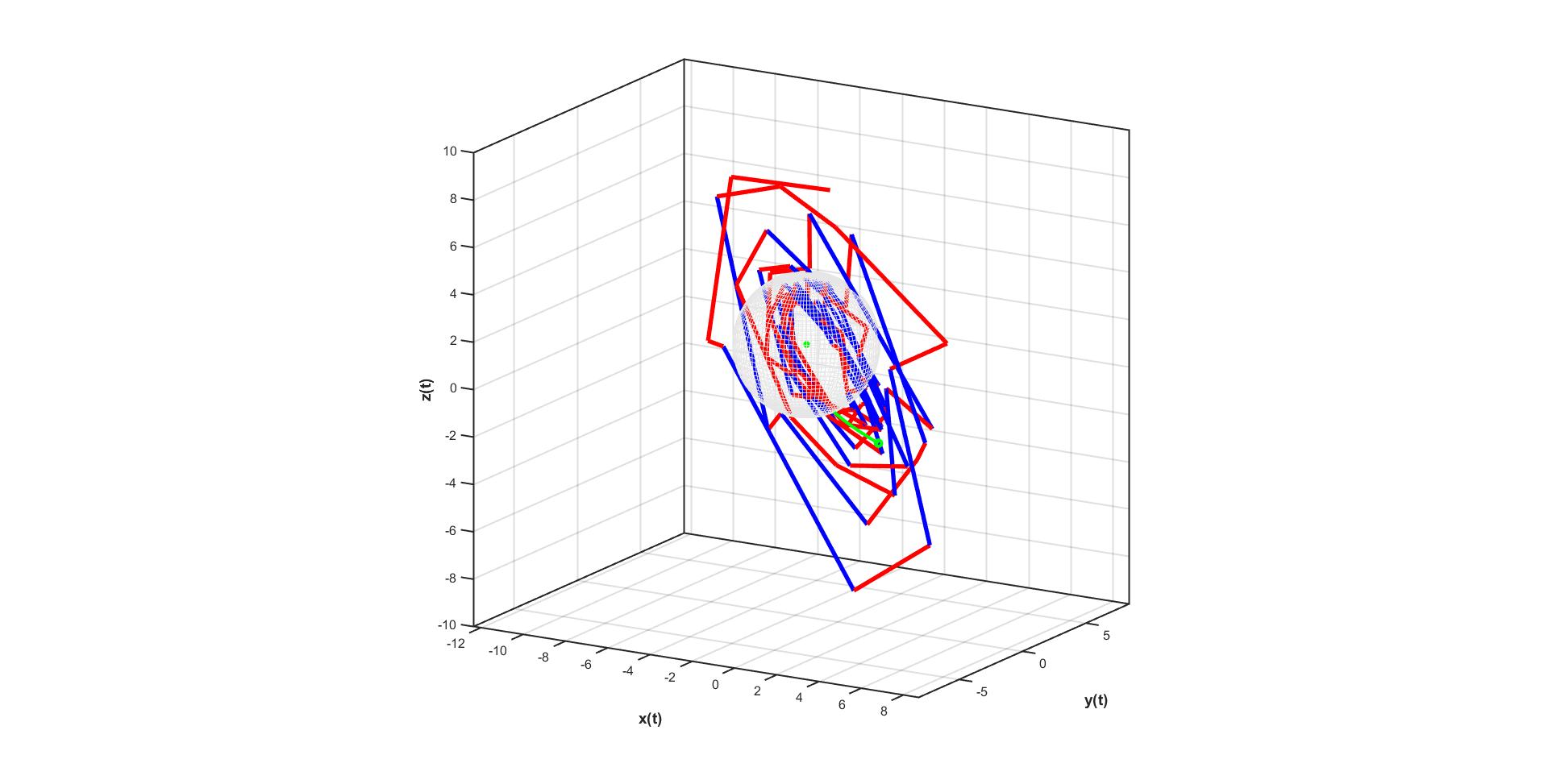}
    \end{minipage}%
    \begin{minipage}{0.5\textwidth}
        \centering
        \includegraphics[width=1\linewidth, height=0.15\textheight]{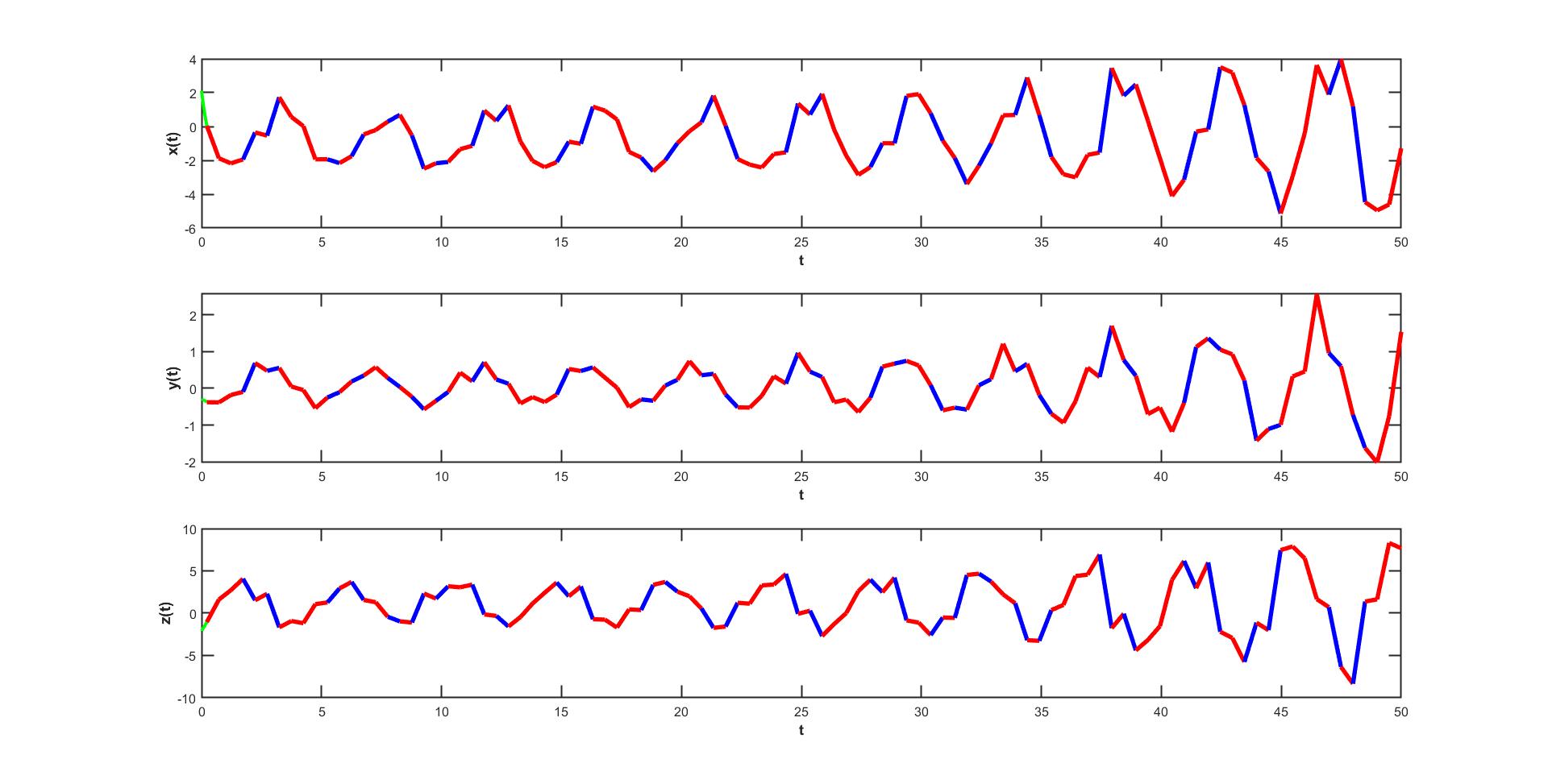}\\
        \includegraphics[width=1\linewidth, height=0.15\textheight]{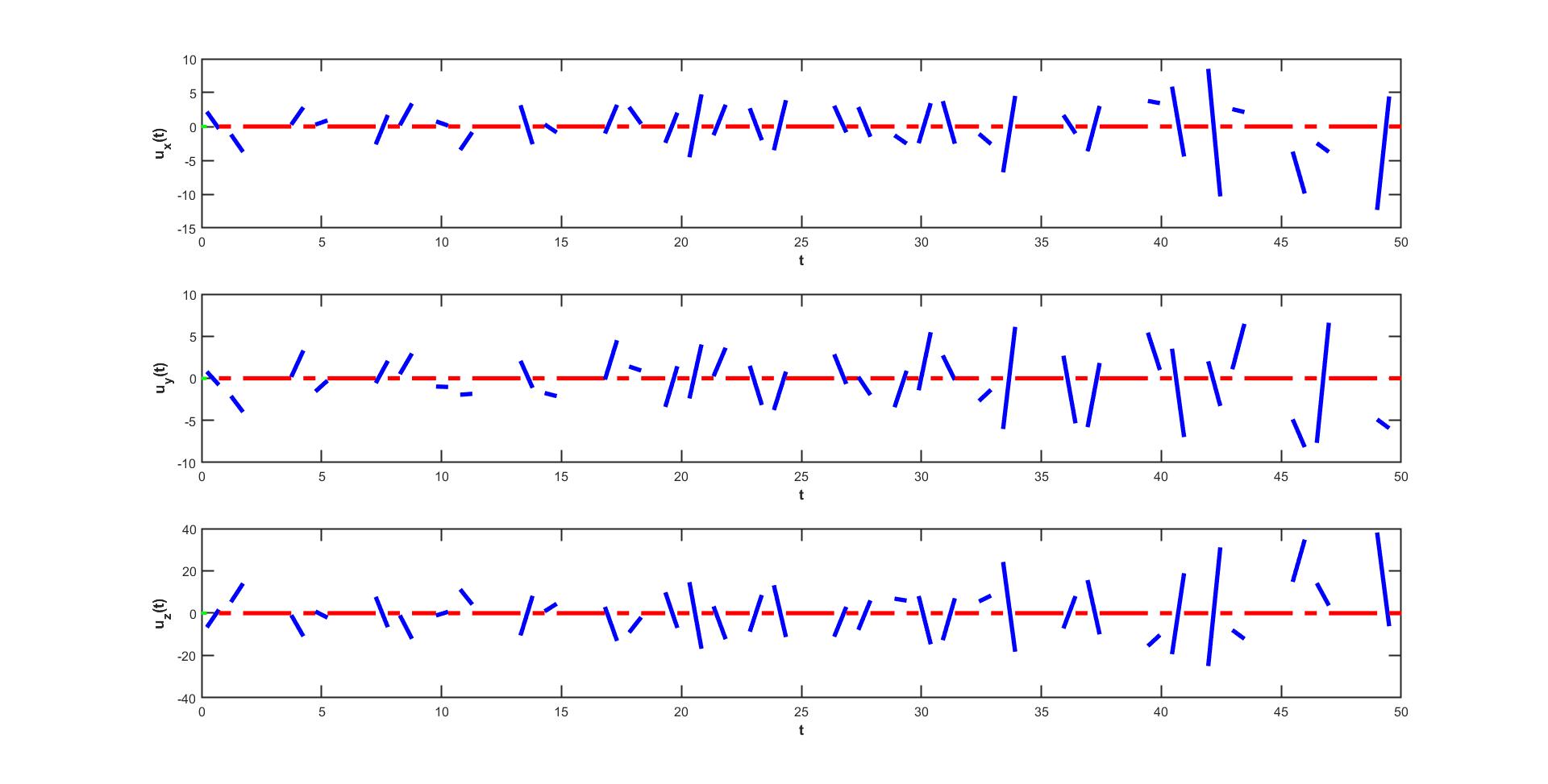}        
    \end{minipage}
    \caption{Double scroll controlled with $\rho=0.6$, $\theta \approx 1.5100$, $\tau = 0.29$ and $\delta=3$. Trajectory initiated in  $CI=(2.104,\ -0.3188,\ -2.0866)$: a) in 3D-space, b) state and control signals.}
    \label{Figure 8}
\end{figure}


\subsection{Rossler system}
 
The well known Rossler system is a nonlinear system determined by the following set of continuously  differentiable equations:

\begin{equation*} \left\{
\begin{array}{l l l} 
\dot{x}&=&-y-z\\
\dot{y}&=&x+0.2y\\
\dot{z}&=&0.2+z(x-5.7)	
\end{array} \right. \end{equation*}

For $a=0.2$, $b=0.2$ and, $c=5.7$, this system is chaotic. This case is considered in \cite{GUAN2006} for stabilizing it in its equilibrium point, by means of a delayed feedback control design.

Here, the objective is to control chaos by stabilizing the system in its equilibrium point $P=(0.007,\ -0.0352,\ 0.0352)$ which linearization matrix has the eigenvalues $\lambda_1 \approx -5.6940$, $\lambda_{2,3} \approx 0.0970\pm 0.9952i$.

For control implementation, it is chosen  $\rho=0.2$, $\theta=\arctan (\omega/\mu) \approx 1.4736$ where $\mu=$ Re$(\lambda_2)$; $\omega=$Im$(\lambda_2)$ and, $\tau=0.2$, and so the matrix $K$  results: 
\begin{equation*}
K= \begin{bmatrix} 
2.9741 & -22.4906 & -1.1580\\
22.4906 & 7.5995 & -3.6157\\
0.0408 & -0.1273 & -0.0106\\ 
     \end{bmatrix} 
\end{equation*}

Putting  $\delta = 5$, the results of the simulation with initial condition $CI=(10,\ 10,\ 10)$ are displayed in Figures \ref{Fig 9} and \ref{Fig 10}.

\begin{figure*} 
    \includegraphics[width=0.55\textwidth]{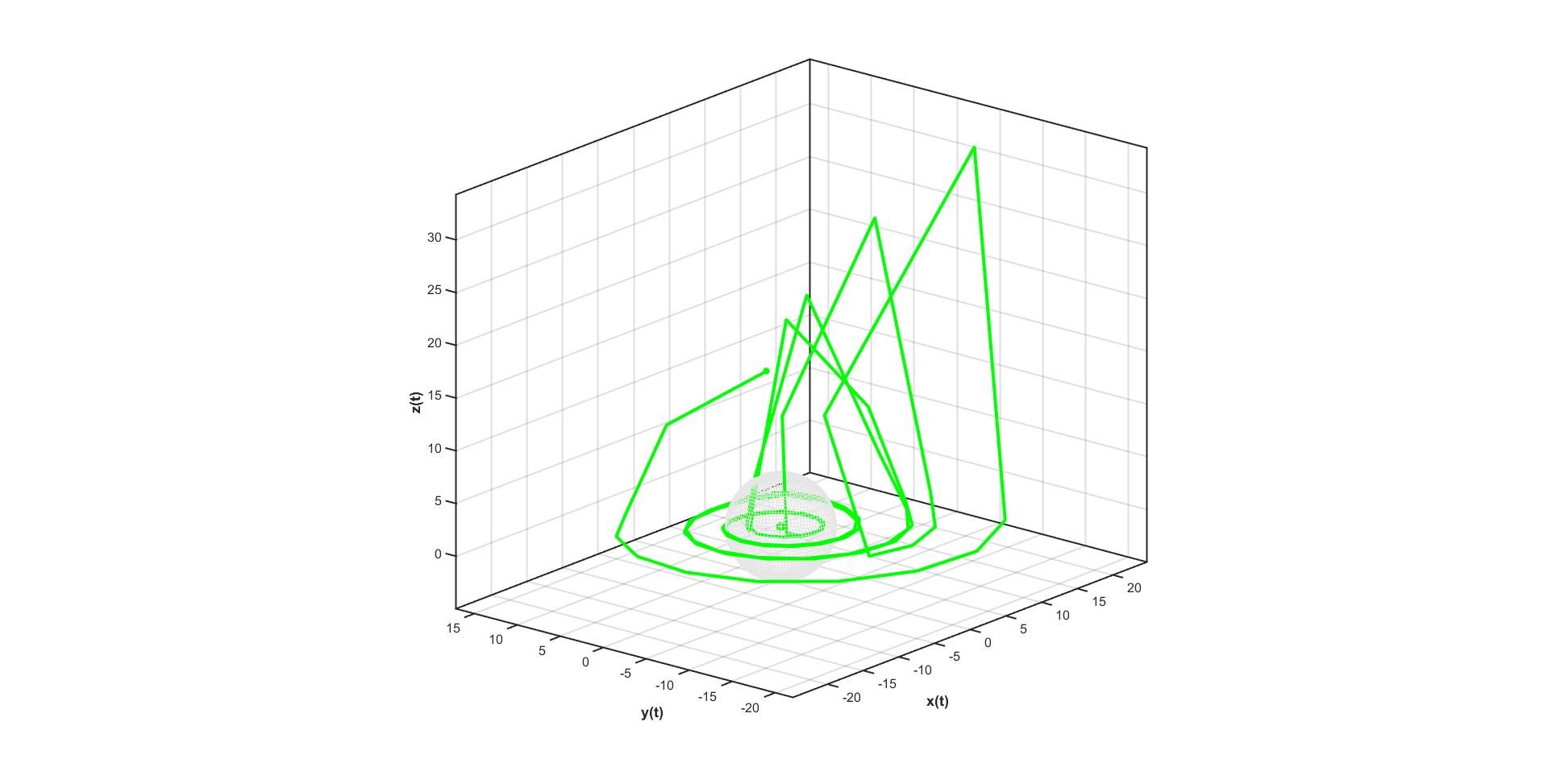}
    \includegraphics[width=0.55\textwidth]{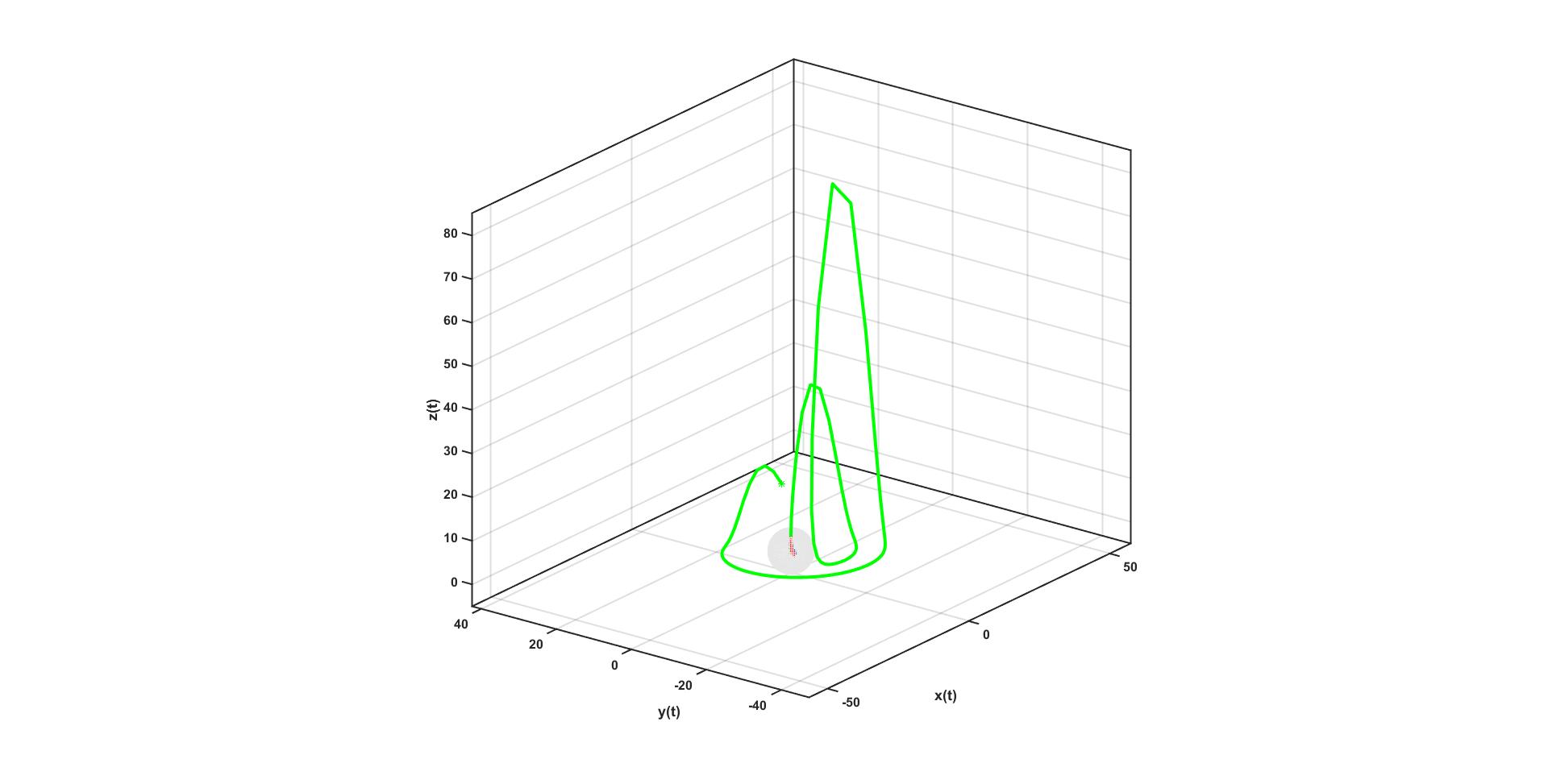}
    \caption{Rossler. Trajectory initialized in CI = (10, 10, 10) of
a) the free and b) the controlled system with $\tau = 0.25, \rho=0.2, \theta \approx 1.4736$ and $\delta = 5$.}
\label{Fig 9}
\end{figure*}
\begin{figure*} 
    \includegraphics[width=0.5\textwidth]{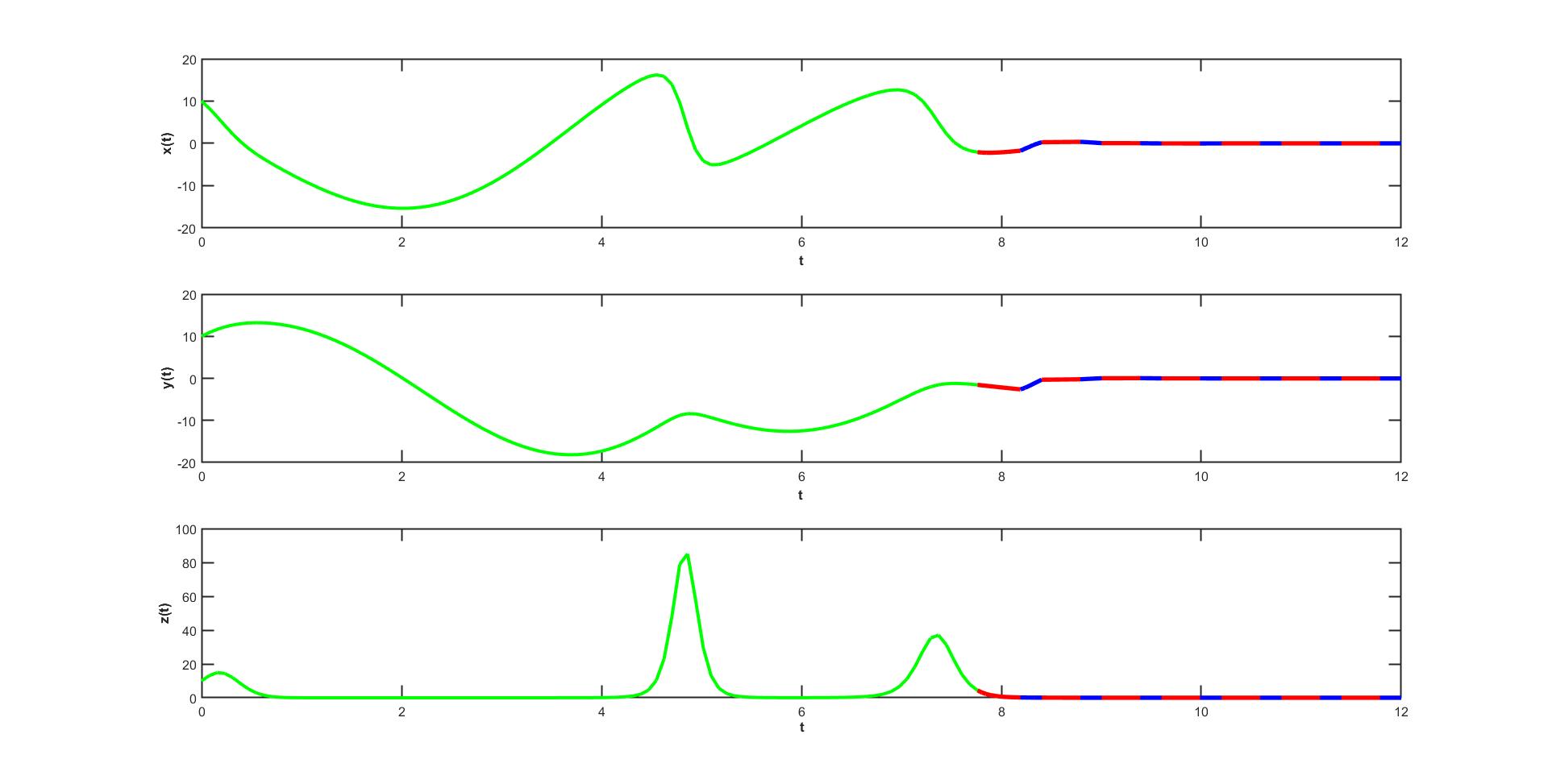}
    \includegraphics[width=0.5\textwidth]{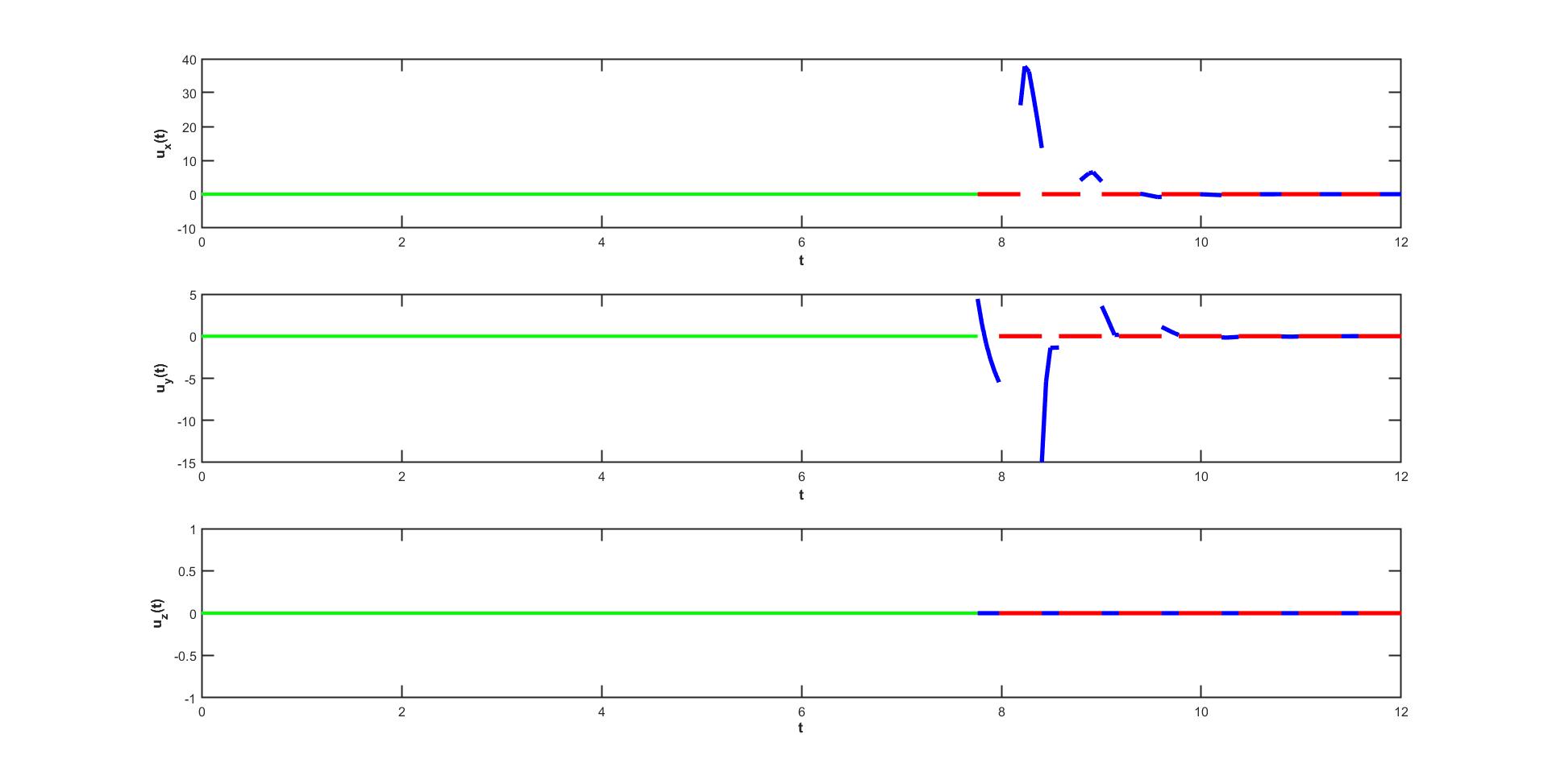}
    \caption{Rossler. State and control signals of the controlled system initialized in CI= (10, 10, 10) with $\tau = 0.25,\ \rho=0.2,\ \theta \approx 1.47363$ and $\delta = 5$. }
\label{Fig 10}
\end{figure*}   

\medskip

Remark: This example is dealt with in \cite{GUAN2006} where the control objective is satisfied by applying a full delayed feedback control. The resulting signals (see Figures \ref{Fig7-DSFE} to \ref{Fig 10} of \cite{GUAN2006} are comparable to the ones obtained by this work proposal (Figures \ref{Fig 9} and \ref{Fig 10}). In \cite{GUAN2006}, controller gain and time delay are obtained by solving certain linear matrix inequalities; moreover, the time delay choosing restricts to an upper bound. On the contrary, for the strategy introduced here, given a desired convergence index, the control gain is computed explicitly for any time delay value.

\section{Concluding remarks}

A control scheme for equilibrium points stabilization in the n-dimensional case has been proposed. It is a full state feedback control based on the difference between  the $\tau$ and $2\tau$-delayed states and an adequate periodic control gain which depends on a constant matrix $K$.

The controlled system results a discontinuous time-delayed differential equation. Its stabilization in a desired equilibrium point is rigorously studied through an associated continuously  differentiable map. It follows that given any hyperbolic equilibrium point:

\begin{itemize}

\item  fixed any pair of stabilizing control parameters $\tau$ and $K$, the control is achieved; in particular, the ONL is overcome, 

\item the set of the stabilizing control parameters is analytically obtained,

\item the control parameters can be systematically computed and their values impact  on different performance aspects. 

\end{itemize}

By adding a nearness-criteria with a third control parameter, $\delta$, a strategy is implemented for controlling chaos by stabilizing equilibrium points embedded in an strange attractor. Performance features of this strategy when applied to the well known chaotic Chua and Rossler systems are observed and some improvements over previous schemes are pointed out.

\section{Future work}

Motivated by the afermentioned advantages of the TDFC with periodic gain, based on full state, novel proposals may be designed with the aim of improving control features in different scenarios:
\begin{itemize}
\item The issue of equilibrium point stabilization for a   control system given by
\begin{equation} \label{eqNLS}
    \dot{x}=f(x,u), \ \ \ y= g(x)
\end{equation} 
where $x \in \mathbb{R}^n $ is the state vector, $u \in \mathbb{R}^m$ is the input (control) vector, $y \in \mathbb{R}^l$ is the output vector, and the functions $f$ and $g$ are differentiable,
has been dealt with by means of time-invariant static output feedback (\cite{Leo-Shuma2010}), by nonstationary static output feedback - also called the Brockett problem - (\cite{Leo-Shuma2010}, \cite{Shumafov}) and by delayed output feedback (\cite{Leo-etal2016}). 

System (\ref{eqNLS}) with null control has an unstable hyperbolic equilibrium point $x^*$, that is $f(x^*,0)=0$. 
The objective is to stabilize the system in $x^*$ by applying an output feedback control.

The delay feedback stabilization for two and three dimensional controllable linear systems is studied in \cite{Leo-etal2016} considering scalar input and scalar output (i.e., $m=l=1$) where  one of the proposed method is based on Pyragas delayed method, that is, 
 $$u(t)= k\big(y(t-\tau)-y(t)\big)$$ 
for adequate $k\neq 0$ and $\tau>0$. Necessary and sufficient conditions on system coefficients for stabilization are obtained by analytical tools in \cite{Leo-etal2016}. Although improvement on the stability regions respect to the static feedback controls is reported, there are cases for which  this methods does not work, namely, the ONL remains. Hence, it sounds attractive to propose the delayed output feedback: 
$$u(t)=\epsilon (t)\big(y(t-2\tau)-y(t-\tau)\big)$$
with $\epsilon$ defined as in (\ref{eq7bis}).

\item In the MDFC proposal (\cite{Ahlborn-Parlitz2005}, \cite{Lu-Li2009}) two (or more delays) take part into the feedback control. It is required that the delays are not integer multiples of each other. By this method, equilibrium point stabilization is achieved for a broader region of stability parameters in comparison to DFC and EDFC. Inspired in this proposal, a suitable modification on the TDFC could be: 
$$u(t)= k(t)\big(y(t-\tau_1)-y(t-\tau_2)\big)$$ 
being $\tau_1<\tau_2$ mutually prime delays  and $k(t)$ an adequate periodic control gain.

\item The design of a strategy as (\ref{eq3}) for UPO´s stabilization  is quite straight. Suppose that $\widetilde{x}(t)$ is a UPO which is not precisely located but its period $T$ is known. The novel oscillating feedback control based on delayed states given by 
 \begin{equation*}
         u(\!t)\!=\! K(\!t)[x(\!t\!-\!2T)\!-\!x(\!t\!-\!T)]
\end{equation*}
with $K(\!t\!)$ an oscillating control gain could be proposed.
By introducing $\delta x=x-\widetilde{x}(t)$, the problem is reformulated as the problem of stabilizing the origin in the non-autonomous n-dimensional case. As the methods developed in \cite{PYRAGAS2018574} and \cite{LEONOV201434}, this proposal involves a periodic control gain but, it differs from them in the difference between delayed states. It is worth to study its achievements and to confront it with the previous schemes. In turn, an  EDFC (\cite{Socolar}) based on the TDFC might be built  for dealing with high-instabilities scenarios. These problems, and additionally, their application for controlling chaos,
by stabilizing UPO's embedded in a strange attractor are interesting for future research.

\end{itemize}

    \nocite{*}
    \bibliography{Ourbiblio}
    \vspace{1.5cm}
\end{document}